# Multiphysics model reduction of thermomechanical vibration in a state-space formulation

Jun-Geol Ahn[1], Hyun-Ik Yang[1*], Jin-Gyun Kim[2*]

[1] *Department of Mechanical Design Engineering, Hanyang University, Wangsibri-ro 222, Seongdong-gu, Seoul, Republic of Korea*
[2] *Department of Mechanical Engineering, Kyung Hee University, 1732 Deogyeong-daero, Giheung-gu, Yongin-si 17104, Republic of Korea*

**Abstract**

The aim of this work is to propose a new multiphysics mode synthesis (MMS) for the thermomechanical vibration problem. The present thermomechanical model is based on a state-space formulation, which consists of displacement $u$, velocity $\dot{u}$, and temperature shift $\theta$. The state-space based thermomechanical formulation is symmetric unlike a conventional non-symmetric $(u, \theta)$ formulation. In the proposed MMS, the structural variables, $u$ and $\dot{u}$, are first reduced, which is then applied to the coupling term in the thermal parts. A term of the thermal domain is then reduced while preserving the multiphysics coupling effects, resulting in improved accuracy. The proposed two-step MMS with the thermal physics domain update can be implemented with the coupling term derived by using the residual flexibility. The proposed MMS strategy can be also applied to accelerate the computational speed by using independent parallel solvers. The performance of the proposed MMS method is evaluated through numerical examples.

[*] Corresponding authors: Hyun-Ik Yang , Jin-Gyun Kim (jingyun.kim@khu.ac.kr)

# 1. Introduction

Considering the thermal and elastic structural coupling effect has been a main issue in many engineering practices [1-3]. In particular, the fully coupled thermomechanical formulation, which can transfer coupling effects from a structural part to a thermal part and vice versa, provides more realistic computational model [4-7]. In the fully coupled thermomechanical problem, the structural and thermal parts are commonly discretized using the finite element (FE) method with the assumptions of thermoelastic state. Among various formulations, a $(u, \theta)$ form with the second order ordinary differential equation (ODE) is perhaps the most widely used wherein $(u, \theta)$ are the structural displacement and thermal shift variables, respectively [4,6,8-9]. While the $(u, \theta)$ formulation performs well for moderate size problems, computational issues of large-scale multiphysics problems have been encountered for design and optimization as well as experimental validation.

Modal projection using a mathematical basis has been a popular remedy to tackle this issue effectively. The coupled eigenvectors obtained from the $(u, \theta)$ formulation is the easiest way to compute the mathematical basis. However, it may not be recommended due to large size of the coupled eigenvalue problems with asymmetric and the resulting coupled modes that cannot be explained in physics. Therefore, normally, structural and thermal eigenvectors are independently carried out from their own single-field eigenvalue problems, whose a few generalized coordinates are then used to reduce respectively the model size of resulting structural and thermal domains. Often, the so-called constraint modes that describe the interaction between two physics are used for accuracy improvement of the resulting reduced models. Such approaches are known as uncoupled and weakly coupled model reduction techniques [10-12]. Krylov subspace and/or singular value decomposition (SVD) may be also considered to compute the basis instead of the eigenvalue problem [12].

It should be noted that the most conventional $(u, \theta)$ formulations are asymmetric, which are less attractive than symmetric formulations in the computational point of view as one must invoke asymmetric linear and/or eigenvalue solvers. The non-symmetric dynamics problems may also cause numerical instability of standard time integration algorithms due to indeterminable eigenvalues arising from interface compatibility conditions. In fact, asymmetric formulations are common in the multiphysics problems [13], and thus it motivated a lot of studies to derive symmetric formulations [14-15]. In the thermomechanical problem, using a space-state approach, symmetric ($u, \dot{u}, \theta$) formulations with the first order ODE can be obtained from the asymmetric ($u, \theta$) formulations [16]. However, the symmetric ($u, \dot{u}, \theta$) formulations that additionally needs the structural velocity ($\dot{u}$) come with larger degrees of freedom (DOFs) than the asymmetric ($u, \theta$) formulations, which requires twice larger DOFs in the structural domain. Using coupled eigenvector is also available to construct reduced model of the symmetric ($u, \dot{u}, \theta$) formulations, but the computational cost of its coupled eigenvalue problem becomes more severe than the ($u, \theta$) formulations. Steenhoek et al. employed the uncoupled model reduction technique, which is based on uncoupled independent structural and thermal eigenvectors, to the ($u, \dot{u}, \theta$) formulation [16]. However, unlike the ($u, \theta$) formulations, common constraint modes containing structural-thermal coupling effects are difficult to apply to the ($u, \dot{u}, \theta$) formulation because the constraint modes derived from Schur complement under the static assumption may not be defined in the state-space formulation with mixed variables [10-12]. It may cause less accurate coupled modes resulting the uncoupled model reduction of the ($u, \dot{u}, \theta$) formulation.

We here propose a two-step multiphysics mode synthesis (MMS) of the symmetric ($u, \dot{u}, \theta$) formulation with a new coupling strategy to handle the issue. The first reduction of the proposed MMS starts from the structural variables, $u$ and $\dot{u}$, by using a few structural modes, which are computed by solving same structural eigenvalue problems to the case of the uncoupled model reduction. Specifically, structural residual flexibility that contains the truncated modal effect of the first reduction is then used as a bridge to transfer the structural modal energy to the thermal domain. As a result, the thermal system matrices are updated by structural modal effect. The reduction of thermal variable $\theta$ using thermal modes is then followed. It should be noted that the thermal modes used in the proposed MMS, which are obtained from the updated thermal eigenvalue problem, are more precisely described features of the thermomechanical problems than the independent thermal modes of the uncoupled model reduction. This is a key role in the proposed MMS with better accuracy. The proposed MMS also possesses an advantage of computational efficiency due to handling symmetric system matrices. In addition, the proposed MMS has a lot of possibilities to accelerate its computational speed by using the single-field parallel solvers since it is derived in the manner of the partitioned method [23-24].

It should be noted that, in a mono physics problem, the residual flexibility is popular in the free interface component mode synthesis (CMS) methods [17-20], and has been also widely used to improve the accuracy of the fixed interface CMS methods [21-22]. However, in the proposed method with the multiphysics problem, the residual flexibility is employed to make a coupling matrix between the structural and thermal domains. Therefore, the usage of the residual flexibility in this work is distinguished from the well-known approaches in the mono-physics problems [17-22], and this work is the first approach, a coupling matrix with the residual flexibility, in the multiphysics model reduction in our knowledge.

The manuscript is organized as follows: a problem setup of the fully coupled thermoelastic problem and its state-space form are reviewed in Section 2. The conventional uncoupled mode reduction method is briefly described in Section 3, and the proposed two-step model reduction method is explained in Section 4. Section 5 shows the performance of the proposed method through numerical examples. Conclusions are given in Section 6.

**2. State-space form of fully coupled thermoelastic problem.**

In this section, we consider a FE discretization form for the fully coupled thermoelastic problem in an object structure. To easily explain the fully coupled thermoelastic problem and its FE discretization, we briefly describe an example structure $\boldsymbol{\Omega}$ as shown in Fig. 1. Here, $A$ and $V$ denote the surface area and volume of the structure $\boldsymbol{\Omega}$, respectively. The unit normal vector of the structure is represented using $\mathbf{n}$. The mechanical material properties $\rho$, $E$, and $v$ represent the density, Young's modulus and Poisson's ratio of the structure, respectively. Also, the thermal material properties $c_E$, $\kappa$ and $\alpha$ respectively represent the specific heat per the unit volume, thermal conductivity and thermal expansion coefficient. From the constitutive equations of the thermomechanical problem and the Fourier's law for the heat conduction, the following relationships can be obtained [4]:

$$\sigma_{ij} = C_{ijkl}\varepsilon_{kl} - \beta_{ij}\theta, \tag{1a}$$

$$\zeta = \beta_{ij}\varepsilon_{ij} + \frac{c_E}{T_0}\theta, \tag{1b}$$

$$q_i = -k_{ij}\theta_{,j}, \tag{1c}$$

$$\theta = T - T_0, \tag{1d}$$

where $\sigma_{ij}$ and $\varepsilon_{ij}$ denote the Cauchy stress tensor and linear strain tensor, respectively, and also $\zeta$, $T_0$, $T$ and $q_i$ indicate the entropy per unit volume, an equilibrium temperature, a temperature variable and a heat flux tensor, respectively. Note that the governing equations and their discretization form are expressed using a temperature shift $\theta$ to simplify the expression.

To set up the fully coupled thermoelastic problem, the fourth order tensor $C_{ijkl}$, the second order tensor $\beta_{ij}$ and $k_{ij}$ are assumed as follows [4]:

$$C_{ijkl} = \lambda \delta_{ij}\delta_{kl} + \mu(\delta_{ik}\delta_{jl} + \delta_{il}\delta_{jk}), \tag{2a}$$

$$\beta_{ij} = \alpha(3\lambda + 2\mu)\delta_{ij}, \tag{2b}$$

$$k_{ij} = \kappa \delta_{ij}, \tag{2c}$$

in which $\lambda$ and $\mu$ denote the Lame constants of the object structure. Here, $\delta$ in Eq. (2) denotes the Kronecker delta.

In order to derive the fully coupled thermoelastic problem, the local balance equation of entropy and equation of motion are both considered [4]. The local balance equation of entropy is written as

$$T_0 \dot{\zeta} + q_{i,i} = 0, \tag{3}$$

and the equation of motion is also written as

$$\sigma_{ij,j} + f_i - \rho \ddot{u}_i = 0. \tag{4}$$

Using the integration form and the divergence theorem, the local balance equation of entropy in Eq. (3) can be written as follows:

$$\int_V \delta\theta T_0 \dot{\zeta}\, dV - \int_V \delta\theta_{,i} q_i\, dV + \int_A \delta\theta q_i n_i dA = 0. \tag{5}$$

The same procedure is applied to the equation of motion in Eq. (4) as follows:

$$-\int_V \delta\varepsilon_{ij}\sigma_{ij}\,dV + \int_A \delta u_i \sigma_{ij} n_j\,dA + \int_V \delta u_i(f_i - \rho\ddot{u}_i)dA = 0. \tag{6}$$

Subtrating Eq. (6) from Eq. (5), the following thermoelastic boundary value problem is obtained as follows:

$$\int_V \delta\theta T_0 \dot{\zeta}\,dV - \int_V \delta\theta_{,i} q_i\,dV + \int_V \delta\varepsilon_{ij}\sigma_{ij}\,dV$$
$$= -\int_A \delta\theta q_i n_i dA + \int_A \delta u_i \sigma_{ij} n_j\,dA + \int_V \delta u_i(f_i - \rho\ddot{u}_i)dA. \tag{7}$$

Finally, substituting Eq. (1) into Eq. (7), the standard boundary value problem for the fully coupled thermoelastic structure can be written as follows [4]:

$$T_0 \int_V \delta\theta \beta_{ij}\dot{\varepsilon}_{ij}\,dV + \int_V \delta\theta c_E \dot{\theta}\,dV + \int_V \delta\theta_{,i} k_{ij}\theta_{,j}\,dV + \int_V \delta\varepsilon_{ij} C_{ijkl}\varepsilon_{kl}\,dV$$
$$- \int_V \delta\varepsilon_{ij}\beta_{ij}\theta\,dV + \int_V \delta u_i \rho\ddot{u}_i\,dV = -\int_A \delta\theta Q\,dA + \int_A \delta u_i f_i^S\,dA + \int_V \delta u_i f_i^B\,dV, \tag{8}$$

where $u$ denotes the displacement of the structure. Here $Q$ indicates the thermal excitation, and also $f_i^S$ and $f_i^B$ indicate the surface and body excitations of the structure, respectively.

Using the FE discretization process, the well-known equations of the fully coupled thermoelastic problem can be derived from Eqs. (2) and (8) as follows [4]:

$$\begin{bmatrix} M_{ss} & 0 \\ 0 & 0 \end{bmatrix}\begin{bmatrix} \ddot{u}_s \\ \ddot{\theta}_T \end{bmatrix} + \begin{bmatrix} 0 & 0 \\ D_{Ts} & D_{TT} \end{bmatrix}\begin{bmatrix} \dot{u}_s \\ \dot{\theta}_T \end{bmatrix} + \begin{bmatrix} K_{ss} & -K_{sT} \\ 0 & K_{TT} \end{bmatrix}\begin{bmatrix} u_s \\ \theta_T \end{bmatrix} = \begin{bmatrix} f_s \\ Q_T \end{bmatrix}, \tag{9a}$$

$$\boldsymbol{u}_s \in \mathbb{R}^{N_s \times 1}, \boldsymbol{\theta}_T \in \mathbb{R}^{N_T \times 1}, \tag{9b}$$

and, the component matrices and vectors in Eq. (9) are

$$\int_V \delta u_i \rho \ddot{u}_i\,dV \Rightarrow \delta \boldsymbol{u}_s^T \boldsymbol{M}_{ss} \ddot{\boldsymbol{u}}_s, \tag{10a}$$
$$T_0 \int_V \delta\theta \beta_{ij}\dot{\varepsilon}_{ij}\,dV \Rightarrow \delta \boldsymbol{\theta}_T^T \boldsymbol{D}_{Ts} \dot{\boldsymbol{u}}_s, \tag{10b}$$
$$\int_V \delta\theta c_E \dot{\theta}\,dV \Rightarrow \delta \boldsymbol{\theta}_T^T \boldsymbol{D}_{TT} \dot{\boldsymbol{\theta}}_T, \tag{10c}$$
$$\int_V \delta\varepsilon_{ij} C_{ijkl}\varepsilon_{kl}\,dV \Rightarrow \delta \boldsymbol{u}_s^T \boldsymbol{K}_{ss} \boldsymbol{u}_s, \tag{10d}$$

$$\int_V \delta\varepsilon_{ij}\beta_{ij}\theta \ \mathrm{d}V \Rightarrow \delta\boldsymbol{u}_s^T \boldsymbol{K}_{ST}\boldsymbol{\theta}_T, \tag{10e}$$

$$\int_V \delta\theta_{,i}k_{ij}\theta_{,j} \ \mathrm{d}V \Rightarrow \delta\boldsymbol{\theta}_T^T \boldsymbol{K}_{TT}\boldsymbol{\theta}_T, \tag{10f}$$

$$\int_A \delta u_i f_i^S \mathrm{d}A + \int_V \delta u_i f_i^B \mathrm{d}V \Rightarrow \delta\boldsymbol{u}^T \boldsymbol{f}_s, \tag{10g}$$

$$-\int_A \delta\theta Q \mathrm{d}A \Rightarrow \delta\boldsymbol{\theta}_T^T \boldsymbol{Q}_T, \tag{10h}$$

$$\boldsymbol{M}_{ss} = \boldsymbol{M}_{ss}^T, \ \boldsymbol{K}_{ss} = \boldsymbol{K}_{ss}^T, \ \boldsymbol{D}_{TT} = \boldsymbol{D}_{TT}^T, \ \boldsymbol{K}_{TT} = \boldsymbol{K}_{TT}^T, \tag{10i}$$

in which $\boldsymbol{M}$, $\boldsymbol{D}$ and $\boldsymbol{K}$ respectively represent the mass, damping and stiffness matrices, and an operator $(\cdot)^T$ indicates a transpose of a matrix (or vector). Subscripts $s$ and $T$ denote the structural and thermal part, respectively. The structural and thermal excitation vectors are specified using $\boldsymbol{f}_s$ and $\boldsymbol{Q}_T$, respectively. $\boldsymbol{u}_s$ and $\boldsymbol{\theta}_T$ denote the structural displacement vector and thermal temperature shift vector, respectively. The detailed derivation process from Eqs. (1) to (10) can be found in Ref [4].

The fully coupled thermomechanical formulation in Eq. (9) can be explained by using the local balanced equation of entropy and the Cauchy stress [4, 16, 34]. The local balance equation of entropy in Eq. (3) contains the first order time differentiation for an entropy per volume $\zeta$. Thus, the structural effects can be transmitted to the thermal part using the term of the first order time differentiation, the coupling matrix $\boldsymbol{D}_{Ts}$, see Eq. (9). On the other hand, the equation of motion in Eq. (4) contains the Cauchy stress $\sigma_{ij}$, and then thermal effects can be transmitted to the structural part through the coupling matrix $\boldsymbol{K}_{sT}$ in Eq. (9).

The conventional fully coupled thermoelastic formulation in Eq. (9) is asymmetric, and then the state-space form is an easiest way to make it to the symmetric form [15, 25-26]. The state-space form of Eq. (9) is written as [16]:

$$\boldsymbol{A}\dot{\boldsymbol{d}} + \boldsymbol{B}\boldsymbol{d} = \boldsymbol{f}, \tag{11a}$$

$$\boldsymbol{A} = \begin{bmatrix} -\boldsymbol{K}_{ss} & 0 & 0 \\ 0 & \boldsymbol{M}_{ss} & 0 \\ 0 & 0 & \boldsymbol{D}_{TT} \end{bmatrix}, \ \boldsymbol{B} = \begin{bmatrix} 0 & \boldsymbol{K}_{ss} & 0 \\ \boldsymbol{K}_{ss} & 0 & -\boldsymbol{K}_{sT} \\ 0 & \boldsymbol{D}_{Ts} & \boldsymbol{K}_{TT} \end{bmatrix}, \tag{11b}$$

$$\boldsymbol{f} = \begin{bmatrix} \boldsymbol{f}_s \\ 0 \\ \boldsymbol{Q}_T \end{bmatrix}, \ \boldsymbol{d} = \begin{bmatrix} \boldsymbol{u}_s \\ \dot{\boldsymbol{u}}_s \\ \boldsymbol{\theta}_T \end{bmatrix}, \tag{11c}$$

$$\boldsymbol{A} = \boldsymbol{A}^T, \tag{11d}$$

$$d \in \mathbb{R}^{(2N_s+N_T) \times 1}. \tag{11e}$$

Here, the vector $d$ denotes the response vector of the state-space model. The system matrix $A$ is symmetric matrix, but the matrix $B$ is asymmetric due to the coupling term $K_{ST}$ and $D_{Ts}$. In order to make $B$ to a symmetric form, the following relationship, which is derived from Eqs. (10b) and (10e), is required:

$$K_{Ts} = \frac{1}{T_0} D_{Ts} = K_{sT}^T. \tag{12}$$

Substituting Eq. (12) into Eq. (11), the following symmetrized state-space form can be obtained [16]:

$$A\dot{d} + Bd = f, \tag{13a}$$

$$A = \begin{bmatrix} -K_{ss} & 0 & 0 \\ 0 & M_{ss} & 0 \\ 0 & 0 & -\widehat{D}_{TT} \end{bmatrix}, \quad B = \begin{bmatrix} 0 & K_{ss} & 0 \\ K_{ss} & 0 & -K_{sT} \\ 0 & -K_{Ts} & -\widehat{K}_{TT} \end{bmatrix}, \quad f = \begin{bmatrix} 0 \\ f_s \\ -\widehat{Q}_T \end{bmatrix}, \quad d = \begin{bmatrix} u_s \\ \dot{u}_s \\ \theta_T \end{bmatrix}, \tag{13b}$$

$$\widehat{D}_{TT} = \frac{1}{T_0} D_{TT}, \quad \widehat{K}_{TT} = \frac{1}{T_0} K_{TT}, \quad \widehat{Q}_T = \frac{1}{T_0} Q_T, \quad A = A^T, \quad B = B^T. \tag{13c}$$

In this paper, the symmetric form of Eq. (13) is called the full model for convenience reading.

### 3. Review of the conventional uncoupled model reduction method

The model reduction process is required for effective computation of the fully coupled thermoelastic problem. It becomes more important in the state-space form because its DOFs are larger than the conventional asymmetric form in Eq. (9). In this section, we briefly review the conventional model reduction method, namely the uncoupled method [16].

In the uncoupled reduction method, the displacement vector $u_s$ and temperature shift vector $\theta_T$ are in Eq. (13) are written as

$$u_s = \Phi_s q_s = \Phi_s^d q_s^d + \Phi_s^e q_s^e, \tag{14a}$$

$$\theta_T = \Xi_T r_T = \Xi_T^d r_T^d + \Xi_T^e r_T^e, \tag{14b}$$

$$q_s^d \in \mathbb{R}^{N_s \times N_s^d}, \ q_s^e \in \mathbb{R}^{N_s \times N_s^e}, \ r_T^d \in \mathbb{R}^{N_T \times N_T^d}, \ r_T^e \in \mathbb{R}^{N_T \times N_T^e} \tag{14c}$$

where $\boldsymbol{\Phi}_s$ and $\boldsymbol{\Xi}_T$ respectively denote the eigenvector matrices of the structural and thermal parts, which are independently computed from their own eigenvalue problems. $\boldsymbol{q}_s$ and $\boldsymbol{r}_T$ are the generalized coordinate vectors associated with $\boldsymbol{\Phi}_s$ and $\boldsymbol{\Xi}_T$, respectively. The superscript $d$ denotes the dominant (or selected terms), and the superscript $e$ denotes the residual term [13,18,27].

In the uncoupled method, the residual terms are eliminated, and thus the displacement vector $\boldsymbol{u}_s$ and the temperature shift vector $\boldsymbol{\theta}_T$ are approximated as follows:

$$\boldsymbol{u}_s \approx \boldsymbol{\Phi}_s^d \boldsymbol{q}_s^d, \ \boldsymbol{\theta}_T \approx \boldsymbol{\Xi}_T^d \boldsymbol{r}_T^d. \tag{15}$$

Here, the eigenvector matrices $\boldsymbol{\Phi}_s^d$ and $\boldsymbol{\Xi}_T^d$ can be obtained by solving the generalized eigenvalue problems of their own physical domains as follows:

$$\boldsymbol{K}_{ss}(\boldsymbol{\varphi}_s)_i = (\lambda_s)_i \boldsymbol{M}_{ss}(\boldsymbol{\varphi}_s)_i, \ i = 1 \cdots N_s^d, \tag{16a}$$

$$\widehat{\boldsymbol{K}}_{TT}(\boldsymbol{\xi}_T)_j = (\gamma_T)_j \widehat{\boldsymbol{D}}_{TT}(\boldsymbol{\xi}_T)_j, \ j = 1 \cdots N_T^d, \tag{16b}$$

$$(\boldsymbol{\varphi}_s)_i \in \mathbb{R}^{N_s \times 1}, \ (\lambda_s)_i \in \mathbb{R}^{1 \times 1}, \ (\boldsymbol{\xi}_T)_j \in \mathbb{R}^{N_T \times 1}, \ (\gamma_T)_j \in \mathbb{R}^{1 \times 1}, \tag{16c}$$

in which $((\lambda_s)_i, (\boldsymbol{\varphi}_s)_i)$ and $((\gamma_T)_j, (\boldsymbol{\xi}_T)_j)$ respectively represent the $i$-th eigenpair of the structural part and $j$-th eigenpair of the thermal part. Also, $N_s^d$ and $N_T^d$ denote the numbers of the dominant structural and thermal modes, respectively. In general, the number of $N_s^d$ and $N_T^d$ are much smaller than $N_s$ and $N_T$, respectively, i.e. $N_s^d \ll N_s$ and $N_T^d \ll N_T$. As the results, $\boldsymbol{\Phi}_s^d$ contains $(\boldsymbol{\varphi}_s)_i$ vectors, and thus the matrix size is $N_s \times N_s^d$. Also, $\boldsymbol{\Xi}_T^d$ contains $(\boldsymbol{\xi}_T)_j$ vectors, and thus the matrix size is $N_T \times N_T^d$.

Using the approximation relationships in Eq. (15), the response vector $\boldsymbol{d}$ in Eq. (13) can be directly approximated as follows:

$$\boldsymbol{d} = \begin{bmatrix} \boldsymbol{u}_s \\ \dot{\boldsymbol{u}}_s \\ \boldsymbol{\theta}_T \end{bmatrix} \approx \boldsymbol{T}_U \boldsymbol{d}_U, \tag{17a}$$

$$T_U = \begin{bmatrix} \Phi_S^d & 0 & 0 \\ 0 & \Phi_S^d & 0 \\ 0 & 0 & \Xi_T^d \end{bmatrix}, d_U = \begin{bmatrix} q_S^d \\ \dot{q}_S^d \\ r_T^d \end{bmatrix}, \tag{17b}$$

$$d_U \in \mathbb{R}^{(2N_S^d + N_T^d) \times 1}, \tag{17c}$$

where the subscript $U$ denotes an uncoupled method, and also $T_U$ represents the transformation matrix of the uncoupled method for the state-space model. Using $T_U$ in Eq. (17), the reduced state-space model can be obtained as follows: $f_U = T_U^T f$

$$A_U \dot{d}_U + B_U d_U = f_U, \tag{18a}$$

$$A_U = \begin{bmatrix} -\Lambda_S^d & 0 & 0 \\ 0 & I_S^d & 0 \\ 0 & 0 & -I_T^d \end{bmatrix}, B_U = \begin{bmatrix} 0 & \Lambda_S^d & 0 \\ \Lambda_S^d & 0 & -(\Phi_S^d)^T K_{ST} \Xi_T^d \\ 0 & -(\Xi_T^d)^T K_{TS} \Phi_S^d & -\Theta_T^d \end{bmatrix}, f_U = T_U^T f, \tag{18b}$$

$$\Lambda_S^d = \begin{bmatrix} (\lambda_S^d)_1 & & & 0 \\ & \ddots & & \\ & & (\lambda_S^d)_i & \\ & & & \ddots \\ 0 & & & (\lambda_S^d)_{N_S^d} \end{bmatrix}, i = 1, \cdots, N_S^d, \tag{18c}$$

$$\Theta_T^d = \begin{bmatrix} (\gamma_T^d)_1 & & & 0 \\ & \ddots & & \\ & & (\gamma_T^d)_j & \\ & & & \ddots \\ 0 & & & (\gamma_T^d)_{N_T^d} \end{bmatrix}, j = 1, \cdots, N_T^d, \tag{18d}$$

where $(\lambda_S^d)_i$ and $(\gamma_T^d)_j$ respectively denote the $i$-th dominant eigenvalue for the structural part and $j$-th dominant eigenvalue for the thermal part, respectively. Thus, $(\lambda_S^d)_i$ is associated with the $i$-th column vector of the matrix $\Phi_S^d$, and also $(\gamma_T^d)_j$ is associated with the $j$-th column vector of the matrix $\Xi_T^d$. Diagonal matrices ($\Lambda_S^d$ and $\Theta_T^d$) and the identity matrices ($I_S^d$ and $I_T^d$) in Eq. (18) can be derived from the following orthogonal relationship of the fully coupled thermoelastic problem:

$$(\varphi_s)_i^T K_{ss} (\varphi_s)_j = \begin{cases} (\lambda_s)_i & (i = j) \\ 0 & (i \neq j) \end{cases}, (\varphi_s)_i^T M_{ss} (\varphi_s)_j = \begin{cases} 1 & (i = j) \\ 0 & (i \neq j) \end{cases}, i, j = 1 \cdots N_s, \tag{19a}$$

$$(\xi_T)_i^T \widehat{K}_{TT}(\xi_T)_j = \begin{cases} (\gamma_T)_i (i = j) \\ 0 \ (i \neq j) \end{cases}, \ (\xi_T)_i^T \widehat{D}_{TT}(\xi_T)_j = \begin{cases} 1 \ (i = j) \\ 0 \ (i \neq j) \end{cases}, \ i,j = 1 \cdots N_T. \tag{19b}$$

The sizes of the matrices in Eq. (18) are $(2N_S^d + N_T^d) \times (2N_S^d + N_T^d)$. Consequently, the matrices in Eq. (18) resulting from uncoupled method are successfully reduced to the original state-space form in Eq. (13). However, it should be noted that the off-diagonal terms of $T_U$ in Eq. (17) are all equal to zero, which implies that the structural and thermal variables are approximated without coupling effects. Therefore, the reduced matrices using $T_U$ may lead less accuracy within the thermoelastic problems that structural domain and thermal domain strongly interact with each other. It is a main motivation in this work.

**4. Two-step model reduction with coupling effect**

In this section, we introduce a new reduction method for the fully coupled thermoelastic problem. The proposed method is a two-step reduction process. First, the structural part is reduced using the dominant structural modes, which are used in the uncoupled method. However, unlike the uncoupled method, the residual structural modes that are eliminated in the uncoupled method are considered to update the thermal part. Through this process, the residual structural modes and coupling matrices can be combined into the thermal domain. Therefore, the updated thermal eigenvector matrix can consider the structural model coupling effects as well as the original thermal effects. The thermal part is then reduced by using the dominant thermal modes.

*4.1. Structural part reduction*

The first step of reduction starts from the structural domain. In the proposed method, the residual term of the structural part is not immediately eliminated. Thus, displacement vector $u_s$ in the proposed method can be written as follows:

$$u_s = \Phi_s^d q_s^d + \Phi_s^e q_s^e, \tag{20a}$$

$$q_s^d \in \mathbb{R}^{N_s \times N_s^d}, \ q_s^e \in \mathbb{R}^{(N_s - N_s^d) \times 1}. \tag{20b}$$

Using the whole eigenvectors of Eq. (20) without neglecting residual modes, the response vector $d$ in Eq. (13) is rewritten as

$$d = T_1 d_1, \tag{21a}$$

$$T_1 = \begin{bmatrix} \boldsymbol{\Phi}_S^d & 0 & \boldsymbol{\Phi}_S^e & 0 & 0 \\ 0 & \boldsymbol{\Phi}_S^d & 0 & \boldsymbol{\Phi}_S^e & 0 \\ 0 & 0 & 0 & 0 & I \end{bmatrix}, d_1 = \begin{bmatrix} q_S^d \\ \dot{q}_S^d \\ q_S^e \\ \dot{q}_S^e \\ \boldsymbol{\theta}_T \end{bmatrix}. \tag{21b}$$

Considering harmonic response $\frac{d}{dt} = \hat{\imath}\omega$ within the free vibration ($f = 0$) from Eq. (13), and then applying the transformation matrix $T_1$ in Eq. (21) with the orthogonal relationship in Eq. (19), we have

$$\tilde{A}\hat{\imath}\omega d_1 + \tilde{B} d_1 = 0, \tag{22a}$$

$$\tilde{A} = \begin{bmatrix} -\Lambda_S^d & 0 & 0 & 0 & 0 \\ 0 & I_S^d & 0 & 0 & 0 \\ 0 & 0 & -\Lambda_S^e & 0 & 0 \\ 0 & 0 & 0 & I_S^e & 0 \\ 0 & 0 & 0 & 0 & -\hat{D}_{TT} \end{bmatrix}, \tilde{B} = \begin{bmatrix} 0 & \Lambda_S^d & 0 & 0 & 0 \\ \Lambda_S^d & 0 & 0 & 0 & -(\boldsymbol{\Phi}_S^d)^T K_{ST} \\ 0 & 0 & 0 & \Lambda_S^e & 0 \\ 0 & 0 & \Lambda_S^e & 0 & -(\boldsymbol{\Phi}_S^e)^T K_{ST} \\ 0 & -K_{TS}\boldsymbol{\Phi}_S^d & 0 & -K_{TS}\boldsymbol{\Phi}_S^e & -\hat{K}_{TT} \end{bmatrix}, \tag{22b}$$

in which $(\lambda_S^e)_j$ denotes the $j$-th residual eigenvalue for the structural part. Thus $(\lambda_S^e)_j$ is associated with the $j$-th column vector of the matrix $\boldsymbol{\Phi}_S^e$. Here, $\hat{\imath}$ in Eq. (22) denotes the imaginary unit. Similar to the case of the uncoupled method, the diagonal matrices ($\Lambda_S^d$, $\Lambda_S^e$) and the identity matrices ($I_S^d$, $I_S^e$) in Eq. (22) can be derived from the orthogonal relationship of the fully coupled thermoelastic problem as shown in Eq. (19). Thus, similar to the diagonal matrix $\Lambda_S^d$ in Eq. (18), the diagonal matrix $\Lambda_S^e$ in Eq. (22) is $\Lambda_S^e = \text{diag}((\lambda_S^e)_1 \cdots (\lambda_S^e)_j \cdots (\lambda_S^e)_{N_S^e})$, $j = 1, \cdots, N_S^e$.

Using the 3rd and 4th rows in Eq. (22b), the residual terms of the structural part can be expressed as

$$\begin{bmatrix} q_s^e \\ \dot{q}_s^e \end{bmatrix} = \left( \begin{bmatrix} 0 & \Lambda_s^e \\ \Lambda_s^e & 0 \end{bmatrix} + \hat{\imath}\omega \begin{bmatrix} -\Lambda_s^e & 0 \\ 0 & I_s^e \end{bmatrix} \right)^{-1} \begin{bmatrix} 0 \\ (\Phi_s^e)^T K_{sT} \end{bmatrix} \theta_T. \qquad (23)$$

Using the Taylor expansion, an inverse term in Eq. (23) can be expanded as

$$\left( \begin{bmatrix} 0 & \Lambda_s^e \\ \Lambda_s^e & 0 \end{bmatrix} + \hat{\imath}\omega \begin{bmatrix} -\Lambda_s^e & 0 \\ 0 & I_s^e \end{bmatrix} \right)^{-1} = \begin{bmatrix} 0 & \Lambda_s^e \\ \Lambda_s^e & 0 \end{bmatrix}^{-1} - \hat{\imath}\omega \begin{bmatrix} 0 & \Lambda_s^e \\ \Lambda_s^e & 0 \end{bmatrix}^{-1} \begin{bmatrix} -\Lambda_s^e & 0 \\ 0 & I_s^e \end{bmatrix} \begin{bmatrix} 0 & \Lambda_s^e \\ \Lambda_s^e & 0 \end{bmatrix}^{-1}$$
$$+ O(\omega^2) + \cdots. \qquad (24)$$

Neglecting terms of a higher order than $\omega$, which have a numerically small contribution in the inverse series expansion, Eq. (24) can be approximated as

$$\left( \begin{bmatrix} 0 & \Lambda_s^e \\ \Lambda_s^e & 0 \end{bmatrix} + \hat{\imath}\omega \begin{bmatrix} -\Lambda_s^e & 0 \\ 0 & I_s^e \end{bmatrix} \right)^{-1} \approx \begin{bmatrix} 0 & (\Lambda_s^e)^{-1} \\ (\Lambda_s^e)^{-1} & 0 \end{bmatrix} + \hat{\imath}\omega \begin{bmatrix} -(\Lambda_s^e)^{-2} & 0 \\ 0 & (\Lambda_s^e)^{-1} \end{bmatrix}. \qquad (25)$$

Substituting Eq. (25) into Eq. (23), the residual terms of the structural part can be approximated as follows:

$$\begin{bmatrix} q_s^e \\ \dot{q}_s^e \end{bmatrix} \approx \begin{bmatrix} (\Lambda_s^e)^{-1}(\Phi_s^e)^T K_{sT} \\ 0 \end{bmatrix} \theta_T + \hat{\imath}\omega \begin{bmatrix} 0 \\ (\Lambda_s^e)^{-1}(\Phi_s^e)^T K_{sT} \end{bmatrix} \theta_T. \qquad (26)$$

In the proposed method, the coupling effects during the reduction process can be implemented by using the approximated residual terms of the structural part. To realize the coupling effects, Eq. (26) and the 5th row of Eq. (22b) are combined as follows:

$$-\bar{D}_{TT}\dot{\theta}_T - K_{Ts}\Phi_s^d \dot{q}_s^d - \widehat{K}_{TT}\theta_T = 0, \qquad (27a)$$
$$\bar{D}_{TT} = \widehat{D}_{TT} + \widehat{\Psi}_{TT}, \qquad (27b)$$
$$\widehat{\Psi}_{TT} = K_{Ts}\Phi_s^e(\Lambda_s^e)^{-1}(\Phi_s^e)^T K_{sT}, \quad \widehat{\Psi}_{TT} = \left(\widehat{\Psi}_{TT}\right)^T. \qquad (27c)$$

Here, $\widehat{\Psi}_{TT}$ in Eq. (27) contains $\Phi_s^e(\Lambda_s^e)^{-1}(\Phi_s^e)^T$, which is the unknown residual eigenpair when calculating only the dominant mode. However, $\Phi_s^e(\Lambda_s^e)^{-1}(\Phi_s^e)^T$ in Eq. (27) can be treated using the orthogonal properties of the fully coupled thermoelastic problem. Considering the orthogonal properties in Eq. (19), $\boldsymbol{\Phi}_S^e (\boldsymbol{\Lambda}_S^e)^{-1} (\boldsymbol{\Phi}_S^e)^T$ in Eq. (27) can be easily calculated using the structural stiffness matrix $\boldsymbol{K}_{SS}$ and dominant mode $\boldsymbol{\Phi}_S^d$ as follows [17-22]:

$$\boldsymbol{\Phi}_S^e (\boldsymbol{\Lambda}_S^e)^{-1} (\boldsymbol{\Phi}_S^e)^T = \boldsymbol{K}_{SS}^{-1} - \boldsymbol{\Phi}_S^d (\boldsymbol{\Lambda}_S^d)^{-1} (\boldsymbol{\Phi}_S^d)^T. \tag{28}$$

This is known as the structural residual flexibility, which is a popular way to indirectly consider the residual model effect in the model reduction process. As a result, the reduction result of the structural part in the proposed method can return the following governing equation:

$$\begin{bmatrix} -\boldsymbol{\Lambda}_S^d & 0 & 0 \\ 0 & \boldsymbol{I}_S^d & 0 \\ 0 & 0 & -\bar{\boldsymbol{D}}_{TT} \end{bmatrix} \begin{bmatrix} \dot{\boldsymbol{q}}_S^d \\ \ddot{\boldsymbol{q}}_S^d \\ \dot{\boldsymbol{\theta}}_T \end{bmatrix} + \begin{bmatrix} 0 & \boldsymbol{\Lambda}_S^d & 0 \\ \boldsymbol{\Lambda}_S^d & 0 & -(\boldsymbol{\Phi}_S^d)^T \boldsymbol{K}_{ST} \\ 0 & -\boldsymbol{K}_{TS} \boldsymbol{\Phi}_S^d & -\widehat{\boldsymbol{K}}_{TT} \end{bmatrix} \begin{bmatrix} \boldsymbol{q}_S^d \\ \dot{\boldsymbol{q}}_S^d \\ \boldsymbol{\theta}_T \end{bmatrix} = \begin{bmatrix} 0 \\ 0 \\ 0 \end{bmatrix}. \tag{29}$$

It shows that the thermal domain in Eq. (29) is updated with the structural residual flexibility. The thermal parts are then sequentially reduced by using new thermal eigenpairs obtained from the updated thermal domain.

*4.2. Thermal part reduction*

Unlike the uncoupled method, the thermal matrix $\bar{\boldsymbol{D}}_{TT}$ in Eq. (29) contains structural modal effects, which are realized by $\widehat{\boldsymbol{\Psi}}_{TT}$. Thus, the new eigenvector matrix $\bar{\boldsymbol{\Xi}}_T^d$ that is different to $\boldsymbol{\Xi}_T^d$ computed from Eq. (29) is obtained from the following generalized eigenvalue problem:

$$\widehat{\boldsymbol{K}}_{TT} (\bar{\boldsymbol{\xi}}_T)_j = (\bar{\gamma}_T)_j \bar{\boldsymbol{D}}_{TT} (\bar{\boldsymbol{\xi}}_T)_j, \quad j = 1 \cdots N_T^d, \tag{30a}$$

$$(\bar{\boldsymbol{\xi}}_T)_j \in \mathbb{R}^{N_T \times 1}, \ (\bar{\gamma}_T)_j \in \mathbb{R}^{1 \times 1}. \tag{30b}$$

Using $\bar{\boldsymbol{\Xi}}_T^d$, the reduced temperature shift vector of the proposed method is approximated as follows:

$$\boldsymbol{\theta}_T \approx \bar{\boldsymbol{\Xi}}_T^d \bar{\boldsymbol{r}}_T^d, \ \bar{\boldsymbol{r}}_T^d \in \mathbb{R}^{N_T^d \times 1}. \tag{31}$$

Consequently, the proposed method leads to the reduced equation of motion and the projection matrix $T_P$ as follows:

$$A_P \dot{d}_P + B_P d_P = f_P, \tag{32a}$$

$$A_P = \begin{bmatrix} -\Lambda_S^d & 0 & 0 \\ 0 & I_S^d & 0 \\ 0 & 0 & -I_T^d \end{bmatrix}, \quad B_P = \begin{bmatrix} 0 & \Lambda_S^d & 0 \\ \Lambda_S^d & 0 & -(\Phi_S^d)^T K_{ST} \bar{\Xi}_T^d \\ 0 & -(\bar{\Xi}_T^d)^T K_{TS} \Phi_S^d & -\bar{\Theta}_T^d \end{bmatrix}, \quad f_P = T_P^T f, \tag{32b}$$

$$\bar{\Theta}_T^d = \begin{bmatrix} (\bar{\gamma}_T^d)_1 & & & 0 \\ & \ddots & & \\ & & (\bar{\gamma}_T^d)_j & \\ & & & \ddots \\ 0 & & & (\bar{\gamma}_T^d)_{N_T^d} \end{bmatrix}, \quad j = 1, \cdots, N_T^d, \tag{32c}$$

$$T_p = \begin{bmatrix} \Phi_S^d & 0 & 0 \\ 0 & \Phi_S^d & 0 \\ 0 & 0 & \bar{\Xi}_T^d \end{bmatrix}, \tag{32d}$$

in which $(\bar{\gamma}_T^d)_j$ denotes the $j$-th dominant eigenvalue for the updated thermal part. Thus $(\bar{\gamma}_T^d)_j$ is associated with the $j$-th column vector of the matrix $\bar{\Xi}_T^d$. Also, the subscript $P$ represents the proposed two-step reduction method. Similar to the case of the uncoupled method, the diagonal matrix $\bar{\Theta}_T^d$ and the identity matrix $I_T^d$ in Eq. (32) can be derived from the following orthogonal relationship:

$$(\bar{\xi}_T)_i^T \widehat{K}_{TT} (\bar{\xi}_T)_j = \begin{cases} (\bar{\gamma}_T)_i & (i = j) \\ 0 & (i \neq j) \end{cases}, \quad (\bar{\xi}_T)_i^T \widehat{D}_{TT} (\bar{\xi}_T)_j = \begin{cases} 1 & (i = j) \\ 0 & (i \neq j) \end{cases}, \quad i, j = 1 \cdots N_T. \tag{33}$$

Finally, the response vector $d$ in Eq. (13) is approximated using the projection matrix $T_P$ as follows:

$$d \approx T_P d_P = T_P \begin{bmatrix} q_S^d \\ \dot{q}_S^d \\ \bar{r}_T^d \end{bmatrix}, \quad d_P = \begin{bmatrix} q_S^d \\ \dot{q}_S^d \\ \bar{r}_T^d \end{bmatrix} \in \mathbb{R}^{(2N_S^d + N_T^d) \times 1}. \tag{34}$$

The reduced-state-space forms resulting from the uncoupled and proposed methods (from Eqs. (18) and (32)) are summarized as

$$
\boldsymbol{A}_U = \begin{bmatrix} -\boldsymbol{\Lambda}_S^d & 0 & 0 \\ 0 & \boldsymbol{I}_S^d & 0 \\ 0 & 0 & -\boldsymbol{I}_T^d \end{bmatrix} \rightarrow \boldsymbol{A}_P = \begin{bmatrix} -\boldsymbol{\Lambda}_S^d & 0 & 0 \\ 0 & \boldsymbol{I}_S^d & 0 \\ 0 & 0 & -\boldsymbol{I}_T^d \end{bmatrix} = \boldsymbol{A}_U, \tag{35a}
$$

$$
\boldsymbol{B}_U = \begin{bmatrix} 0 & \boldsymbol{\Lambda}_S^d & 0 \\ \boldsymbol{\Lambda}_S^d & 0 & -(\boldsymbol{\Phi}_S^d)^T \boldsymbol{K}_{ST} \boldsymbol{\Xi}_T^d \\ 0 & -(\boldsymbol{\Xi}_T^d)^T \boldsymbol{K}_{TS} \boldsymbol{\Phi}_S^d & -\boldsymbol{\Theta}_T^d \end{bmatrix} \rightarrow \boldsymbol{B}_P = \begin{bmatrix} 0 & \boldsymbol{\Lambda}_S^d & 0 \\ \boldsymbol{\Lambda}_S^d & 0 & -(\boldsymbol{\Phi}_S^d)^T \boldsymbol{K}_{ST} \bar{\boldsymbol{\Xi}}_T^d \\ 0 & -(\bar{\boldsymbol{\Xi}}_T^d)^T \boldsymbol{K}_{TS} \boldsymbol{\Phi}_S^d & -\bar{\boldsymbol{\Theta}}_T^d \end{bmatrix}. \tag{35b}
$$

It should be noted that, unlike the original thermal eigenvector matrix $\boldsymbol{\Xi}_T^d$ in the uncoupled method, the thermal eigenvector matrix $\bar{\boldsymbol{\Xi}}_T^d$ contains the structural modal effects, see Eqs. (18) and (32). Consequently, this two-step reduction additionally provides modal coupling effects as well as independent modal projection of the two different physics. Thus, one can expect better reduced matrices in this process than the uncoupled method described in Section 3. Nevertheless, the size of the reduced matrices in the proposed method is identical to the one of the uncoupled method. In this work, the formulation is derived by using the state-space form, but the equivalent reduced matrices in Eq. (32) could be also derived from the second order formulation in Eq. (9) as shown in Appendix 1.

*4.3. Computational process*

In Eq. (35), most components of the reduced matrices are diagonal matrices of the structural and thermal eigenvalues and/or identity matrices. It clearly shows that the final reduced form of the proposed method is easy to implement even though the complicated sequential reduction process is used to derive it. In addition, computing the new thermal eigenvector matrix $\bar{\boldsymbol{\Xi}}_T^d$ is only a difference to the uncoupled method, which requires a simple modification of thermal matrix. i.e. $\widehat{\boldsymbol{D}}_{TT} \rightarrow \bar{\boldsymbol{D}}_{TT} = \widehat{\boldsymbol{D}}_{TT} + \widehat{\boldsymbol{\Psi}}_{TT}$. As a result, similar to the uncoupled method, the structural and thermal eigenvectors of the proposed method can be computed by using independent structural and thermal solvers as described in the partitioned method [13,24]. The process of the constructing the reduced model with the proposed method is schematically described in Fig. 2, and also the specific algorithm of the proposed method is explained below.

Algorithm 1. Two-step model reduction
---

Input: Structural matrices ($M_{ss}$ and $K_{ss}$), thermal matrices ($\widehat{D}_{TT}$ and $\widehat{K}_{TT}$), and coupling matrices ($K_{sT}$ and $K_{Ts}$)

Dominant mode numbers of the structural and thermal domains ($N_s^d$ and $N_T^d$)

Output: reduced matrices of the proposed method ($A_p$ and $B_p$)

1. Compute the $N_s^d$ eigenpairs of the structural mode

$$K_{ss}(\varphi_s)_i = (\lambda_s)_i M_{ss}(\varphi_s)_i, \ i = 1, \dots, N_s^d.$$

2. Construct the matrices $\Lambda_s^d$ and $\Phi_s^d$

$$\Lambda_s^d = diag\left((\lambda_s^d)_1, \cdots, (\lambda_s^d)_{N_s^d}\right), \ \Phi_s^d = \left[(\varphi_s^d)_1 \ \cdots \ (\varphi_s^d)_{N_s^d}\right].$$

3. Calculate the updated thermal damping matrix $\bar{D}_{TT}$

$$\bar{D}_{TT} = \widehat{D}_{TT} + \widehat{\Psi}_{TT}, \widehat{\Psi}_{TT} = K_{Ts}\left[K_{ss}^{-1} - \Phi_s^d(\Lambda_s^d)^{-1}(\Phi_s^d)^T\right]K_{sT}.$$

4. Compute the $N_T^d$ eigenvalues of the thermal mode

$$\widehat{K}_{TT}(\bar{\xi}_T)_j = (\bar{\gamma}_T)_j \bar{D}_{TT}(\bar{\xi}_T)_j, \ j = 1 \cdots N_T^d.$$

5. Construct the matrices $\bar{\Theta}_T^d$ and $\bar{\Xi}_T^d$

$$\bar{\Theta}_T^d = diag\left((\bar{\gamma}_T^d)_1, \cdots, (\bar{\gamma}_T^d)_{N_T^d}\right), \ \bar{\Xi}_T^d = \left[(\bar{\xi}_T^d)_1 \ \cdots \ (\bar{\xi}_T^d)_{N_T^d}\right].$$

6. Compute the coupling matrix: $-(\Phi_s^d)^T K_{sT} \bar{\Xi}_T^d$.

7. Construct the reduced system matrices $A_p$ and $B_p$ ▷ Output
---

Note that the thermal matrix $\bar{D}_{TT}$ in the proposed method can be a fully populated matrix unlike the uncoupled formulation with a sparse matrix $\widehat{D}_{TT}$ in Eq. (13). This is a well-known limitation when the residual flexibility is used in the formulation of the model reduction both mono- and multiphysics problems [21-22, 28-29]. For the reason, the proposed formulation requires more computational cost to reduce the thermal part than the uncoupled method. However, due to the updated thermal basis in the proposed method, the accuracy of the reduced matrix can be dramatically improved at the same reduced matrix size, and also computing a few eigenpairs of the thermal domain only needs in the model reduction process. The performances of the proposed method are evaluated through numerical examples in the next section.

## 5. Numerical examples

In this section, we evaluate the performance of the proposed method by comparing it with the uncoupled method. Three numerical examples are covered: a 2d plate, a 3d flange-pipe and a lithium-ion battery in the fully coupled thermoelastic state [30-31]. First, the uncoupled and proposed methods are compared in the frequency domain using the relative eigenvalue error of the state-space model. The reduction methods are also evaluated using the transient responses in the time domain. Here, the Runge-Kutta method, which is a widely used time integrator for a state-space model, is employed to calculate the transient responses [32-33].

The $i$-th eigenvalue of the full model can be calculated using Eq. (13) as follows:

$$\boldsymbol{B}(\boldsymbol{\chi}_F)_i = (\mu_F)_i \boldsymbol{A}(\boldsymbol{\chi}_F)_i, \quad i = 1, \cdots, 2N_s + N_T. \tag{36}$$

where $(\mu)_i$ and $(\boldsymbol{\chi})_i$ respectively represent the $i$-th eigenvalue and eigenvector of the state-space model. The subscript $F$ denotes the full model in Eq. (13).

Similar to the case of the full model, the reduced $i$-th eigenvalue of each reduction methods can also be calculated using Eqs. (18) and (32) as follows:

$$\boldsymbol{B}_U(\tilde{\boldsymbol{\chi}}_U)_i = (\tilde{\mu}_U)_i \boldsymbol{A}_U(\tilde{\boldsymbol{\chi}}_U)_i, \quad i = 1, \cdots, 2N_s + N_T, \tag{37a}$$

$$\boldsymbol{B}_P(\tilde{\boldsymbol{\chi}}_P)_i = (\tilde{\mu}_P)_i \boldsymbol{A}_P(\tilde{\boldsymbol{\chi}}_P)_i, \quad i = 1, \cdots, 2N_s + N_T, \tag{37b}$$

where a notion $(\sim)$ indicates the reduced state.

Using the results in Eqs. (36) and (37), the $i$-th relative eigenvalue error is calculated as follows:

$$\left(e_\mu\right)_i = \frac{|(\mu_F)_i - (\tilde{\mu})_i|}{|(\mu_F)_i|}, \quad i = 1, \cdots, N_\mu, \tag{38}$$

in which $N_\mu$ is the number of calculated eigenvalues. Here, the $i$-th reduced eigenvalue $(\tilde{\mu})_i$ depends on the selected reduction model. For example, the $i$-th relative eigenvalue error for the uncoupled model is calculated from $(\tilde{\mu})_i = (\tilde{\mu}_U)_i$. In the coupled problem in Eq. (36), the purely real eigenvalues can be recognized as the thermal dominant eigenvalues and the other eigenvalues are then the structural dominant eigenvalues. This can be found that the number of purely real eigenvalues is the same as the thermal matrix size $N_T$ in Eq. (16), and thus the eigenvalues in this case are classified into a thermal part in this paper [16,34]. The absolute eigen spectra may be used to define the thermal and structural dominant eigenvalues because those may be clearly distinguished in the macro problem. However, it should be carefully used within the micro problem that the ranges of the thermal and structural absolute eigenvalue spectra can be closely met [16,34]. The separation properties of absolute eigenvalue spectra in the thermomechanical problem can be identified using the highest period of frequency of a structural part and the highest time constant of a thermal part [16, 34]. The computer specification is Inter® Xeon® E5-1620 v4, 3.50GHz, RAM 16.GB. In the numerical examples, the equilibrium temperature $T_0$ is selected as $T_0 = 25$ (°C).

*5.1. 2d plate problem*

Let us consider a 2d plate problem as described in Fig. 3, in which shows the FE discretization with the geometry quantities, boundary and excitation conditions for the thermal and structural parts. A thickness in this problem is selected as 0.001m. Also, the height $h$ and width $l$ are selected as 0.042 m and 0.140m, respectively. The material in the 2d plate problem is selected as the silicon, and its material properties are listed as follows [4]: Young's modulus $E = 162.4\text{GPa}$, Poisson ratio $v = 0.28$, density $\rho = 2330\,\frac{\text{kg}}{m^{-3}}$, thermal expansion coefficient $\alpha = 2.54 \times 10^{-6} K^{-1}$, thermal conductivity $\kappa = 145\,\frac{W}{(m^{-1}K^{-1})}$ and heat capacity $\frac{c_E}{\rho} = \frac{711\,\text{J}}{(K \cdot \text{kg})}$ [4]. In the 2d plate problem, $N_s$ and $N_T$ in Eq. (16) are $N_s = 280$ and $N_T = 140$, respectively. Also, $N_s^d$ and $N_T^d$ in Eq. (16) are selected as $N_s^d = 30$ and $N_T^d = 30$, respectively.

The relative eigenvalue errors between the full and reduced models, which are calculated from Eq. (38), are plotted in Fig. 4. As the results, the proposed method can return more high accuracy than the conventional uncoupled method in the frequency domain.

To obtain the responses of the 2d plate problem, a constant thermal excitation $\boldsymbol{Q}_T = 0.1$ (kW) is applied, and also a sinusoidal structural excitation $\boldsymbol{f}_s = 3 \sin 10t$ (kN) is applied in the $y$- direction with the Runge-Kutta method [32-33]. Maximum temperature change and displacement results are shown in Fig. 5. Figs. 6, 7 and 8 show the temperature and displacement fields of the full model over time, and also show the differences in the temperature and displacement fields between the full and reduced models over time. Those numerical results of the transient analysis clearly demonstrate that the proposed method can return higher accuracy results when compared with the uncoupled method.

*5.2. 3d flange-pipe problem*

Next, a 3d flange-pipe problem is considered. Fig. 9 indicates that the result of the FE discretization and geometry quantities for the 3d flange-pipe problem, and also Fig. 10 shows the boundary and excitation conditions for the thermal and structural parts in the 3d flange-pipe problem. The radius $R_1$ and $R_2$ are $0.23$ m and $0.15$ m, respectively, and the width $l$ is 2.8m. The material properties in the 3d flange-pipe problem are as follows: The material in the 3d flange-pipe problem is selected as the aluminum steel alloy, and its material properties are listed as follows [4]: Young's modulus $E = 162.4$GPa, Poisson ratio $v = 0.28$, density $\rho = 2330 \frac{\text{kg}}{m^{-3}}$, thermal expansion coefficient $\alpha = 2.54 \times 10^{-6} K^{-1}$, thermal conductivity $\kappa = 145 \frac{W}{(m^{-1}K^{-1})}$ and heat capacity $\frac{c_E}{\rho} = \frac{711 \text{ J}}{(K \cdot \text{kg})}$ [4]. In the 3d flange-pipe problem, $N_S$ and $N_T$ in Eq. (16) are $N_S = 6000$ and $N_T = 2000$, respectively. Also, $N_S^d$ and $N_T^d$ in Eq. (16) are selected as $N_S^d = 300$ and $N_T^d = 300$, respectively.

The relative eigenvalue errors between the full and reduced models, which are calculated from Eq. (38), are plotted in Fig. 11. In addition, to obtain the transient responses, a sinusoidal thermal excitation $\boldsymbol{Q}_T = \sin 2.5t$ (kW) is applied, and also a sinusoidal structural excitation $\boldsymbol{f}_s = 10 \sin 4t$ (kN) is applied in the $y$-direction. Maximum temperature change and displacement results are shown in Fig. 12. The temperature field of the full model and differences between the full and reduced models for the temperature field are described in Fig. 13. Similarly, the structural displacement field and differences between the full and reduced models are plotted in Figs. 14, 15 and 16. Those numerical results in the 3d flange pipe problem clearly shows the high accuracy of the proposed method from both the frequency and time domain analyses.

The computational cost to construct reduced models are explicitly compared in Table 1. The mode superposition method, which is a common reduction method in the engineering field, is additionally considered. The table shows that the proposed method needs more computational cost than the uncoupled method, but the proposed method can provide high accurate reduced matrices compared to the uncoupled method as shown in Figs. 11 to 16. The mode superposition method may lead enough accutate model because it is based on the reference eigensolution obtained using the full matrices $A$ and $B$, see Eq. (36). However, the mode superposition method required relatively large computational cost to construct the reduction matrices than the other methods with the independent physics basis. This is a reason that the independent physics basis approaches have been used in the model reduction method of the multiphysics problems [11,13,16].

*5.3. Lithium-ion battery problem*

A commercial pouch lithium-ion battery, which has been widely used in the engineering practices, is implemented to evaluate the reduction methods [30-31]. Fig. 17 briefly outlines the structure of the lithium-ion battery and its mesh results. Materials of the lithium-ion battery are selected as Ref. [30-31]. In the lithium-ion battery problem, $N_S$ and $N_T$ in Eq. (16) are $N_S = 15661$ and $N_T = 5187$, respectively. Also, $N_S^d$ and $N_T^d$ in Eq. (16) are selected as $N_S^d = 850$ and $N_T^d = 850$, respectively. In this numerical example, a structural boundary condition is implemented by fixing the bottom of the container, and a thermal boundary condition is also implemented by maintaining the temperature of the bottom and sides of the container at a temperature 25 (°C). The excitation conditions of the lithium-ion battery problem are realized by Ref [30-31].

The relative eigenvalue errors, which can verify the performance in the frequency domain, between the full and reduced models are plotted in Fig. 18. In the time domain analysis, the Runge-Kutta method is used to obtain the temperature fields and the structural fields [32-33]. The temperature fields, which can describe a temperature distribution of the full model and performances of the reduced models, are plotted in Fig. 19. Figs. 20, 21 and 22 show the structural displacement fields over time, and also show the differences between the full model and the reduced models, which can explain the performance of the reduced models. As the results, it is confirmed that the proposed method can be operated with the excellent performance in the numerical examples.

*5.4. 2d plate problem at the micro-scale*

The thermomechanical problem at the micro-scale is one of the main thermomechanical issues to obtain more realistic responses of the micro-systems [16, 34]. In particular, the thermal and structural absolute eigenvalue spectra may be closely met in the same absolute eigenvalue spectra range at the micro-scale [16, 34]. To design the thermomechanical problem at the micro-scale, we invoke the 2d plate problem in Section 5.1. The material properties are the same as Section 5.1. However, to design a micro-system, a thickness in this problem is selected as 0.1μm. Also, the height $h$ and width $l$ are selected as 4μm and 20μm, respectively.

The absolute eigenvalue spectra of the 2d plate problem for the macro-scale ($h = 0.042m$, $l = 0.140m$) and micro-scale ($h = 4$μm, $l = 20$μm) are plotted in Fig. 23. As can be verified from previous studies, the thermal and structural eigenvalue spectra may be closely met in the same absolute eigenvalue spectra range at the micro-scale [16, 34].

To evaluate the proposed method at the micro-scale, the relative eigenvalue errors are computed using various structural and thermal modes ( $N_s^d = 30, 70, 110$ and $N_T^d = 30, 70, 110$) The results of the relative eigenvalue errors are plotted in Fig. 24. As the results, it is confirmed that the proposed method can also be well-operated at the micro-systems.

**6. Conclusion**

A new reduced-order modeling technique for the fully coupled thermoelastic problem, which is formulated as a state-space form, is proposed and evaluated. The new multiphysics mode synthesis (MMS) is a two-step sequential reduction from the structural to thermal domains unlike the conventional uncoupled method. Structural residual flexibility is then employed to consider the modal coupling effects between the structural and thermal domains, which acts as a canal to transfer the structural modal energy to the thermal domains in the sequential projection process. Numerical results, including a practical lithium-ion battery problem, clearly demonstrate high accuracy of the proposed method from both the frequency and time domain analyses. This work clearly shows that making a transfer canal between different physics and/or different domains is essential to well construct the reduced matrices of multiphysics problems. In this work, we only dealt with the state-space form of the thermoelastic problem, but the main idea of the proposed method can be extended to other multi-physics problems. In particular, the reduced matrices resulting from the proposed method are highly sparse, and also those are simply computed by using independent structural and thermal solvers. Those many imply excellent potentials of the proposed method in the computational aspects by combining high performance parallel computing strategies although it requires an updating process of the thermal domain with relatively condensed matrix computation.

**Acknowledgement**

This research was supported by the Basic Science Research Programs through the National Research Foundation of Korea funded by the Ministry of Science, ICT, and Future Planning (NRF-2018R1A1A1A05078730).

Appendix 1. Two-step model reduction from the second order governing equation.

In this Appendix 1, we derive the proposed method from the second-order form in Eq. (9). Through this process, we want to show that the reduction result in Appendix 1 is equivalent to the proposed formulation in Section 4.

The proposed method using the second-order form is derived from Eqs. (9), (12) and (13) with free vibration as follows:

$$\begin{bmatrix} M_{ss} & 0 \\ 0 & 0 \end{bmatrix} \begin{bmatrix} \ddot{u}_s \\ \ddot{\theta}_T \end{bmatrix} + \begin{bmatrix} 0 & 0 \\ -K_{Ts} & -\widehat{D}_{TT} \end{bmatrix} \begin{bmatrix} \dot{u}_s \\ \dot{\theta}_T \end{bmatrix} + \begin{bmatrix} K_{ss} & -K_{sT} \\ 0 & -\widehat{K}_{TT} \end{bmatrix} \begin{bmatrix} u_s \\ \theta_T \end{bmatrix} = \begin{bmatrix} 0 \\ 0 \end{bmatrix}, \quad K_{Ts} = \frac{1}{T_0} D_{Ts}. \qquad (A1)$$

Similar to the proposed method in Section 4, the displacement vector $u_s$ is written as

$$u_s = \Phi_s^d q_s^d + \Phi_s^e q_s^e. \qquad (A2)$$

Using the orthogonal relationship in Eq. (19), the displacement vector in Eq. (A2) and a first row in Eq. (A1), the following relationships can be obtained:

$$I_s^d \ddot{q}_s^d + \Lambda_s^d q_s^d - (\Phi_s^d)^T K_{sT} \theta_T = 0, \qquad (A3a)$$
$$I_s^e \ddot{q}_s^e + \Lambda_s^e q_s^e - (\Phi_s^e)^T K_{sT} \theta_T = 0. \qquad (A3b)$$

Neglecting a term of higher order in Eq. (A3b), the residual displacement $q_s^e$ can be approximated as

$$q_s^e = (\Lambda_s^e)^{-1} (\Phi_s^e)^T K_{sT} \theta_T. \qquad (A4)$$

Substituting Eq. (A4) into Eq. (A2), and then applying the residual flexibility in Eq. (28), the displacement vector in Eq. (A2) can be approximated as follows:

$$u_s \simeq \Phi_s^d q_s^d + [K_{ss}^{-1} - \Phi_s^d (\Lambda_s^d)^{-1} (\Phi_s^d)^T] K_{sT} \theta_T. \qquad (A5)$$

Substituting Eq. (A5) into second row in Eq. (A1), the thermal part for the proposed method in Eq. (29) can be obtained:

$$-\boldsymbol{K}_{TS}\boldsymbol{\Phi}_S^d \dot{\boldsymbol{q}}_S^d - \bar{\boldsymbol{D}}_{TT}\dot{\boldsymbol{\theta}}_T - \widehat{\boldsymbol{K}}_{TT}\boldsymbol{\theta}_T = \boldsymbol{0}. \tag{A6}$$

Similar to the proposed method, the basis of the thermal part is obtained from the matrices $\widehat{\boldsymbol{K}}_{TT}$ and $\bar{\boldsymbol{D}}_{TT}$, and thus the temperature shift vector $\boldsymbol{\theta}$ can be approximated as in Eq. (31):

$$\boldsymbol{\theta}_T \approx \bar{\boldsymbol{\Xi}}_T^d \bar{\boldsymbol{r}}_T^d, \quad \bar{\boldsymbol{r}}_T^d \in \mathbb{R}^{N_T^d \times 1}. \tag{A7}$$

Therefore, the reduction results of the thermal and structural parts in the second-order form can be obtained using Eqs. (A3), (A6), (A7) and the orthogonal relationship in Eq. (30) as follows:

$$\begin{bmatrix} \boldsymbol{I}_S^d & \boldsymbol{0} \\ \boldsymbol{0} & \boldsymbol{0} \end{bmatrix} \begin{bmatrix} \ddot{\boldsymbol{q}}_S^d \\ \ddot{\boldsymbol{r}}_T^d \end{bmatrix} + \begin{bmatrix} \boldsymbol{0} & \boldsymbol{0} \\ -(\bar{\boldsymbol{\Xi}}_T^d)^T \boldsymbol{K}_{TS}\boldsymbol{\Phi}_S^d & -\boldsymbol{I}_T^d \end{bmatrix} \begin{bmatrix} \dot{\boldsymbol{q}}_S^d \\ \dot{\boldsymbol{r}}_T^d \end{bmatrix} + \begin{bmatrix} \boldsymbol{\Lambda}_S^d & -(\boldsymbol{\Phi}_S^d)^T \boldsymbol{K}_{ST}\bar{\boldsymbol{\Xi}}_T^d \\ \boldsymbol{0} & -\bar{\boldsymbol{\Theta}}_T^d \end{bmatrix} \begin{bmatrix} \boldsymbol{q}_S^d \\ \bar{\boldsymbol{r}}_T^d \end{bmatrix} = \begin{bmatrix} (\boldsymbol{\Phi}_S^d)^T \boldsymbol{f}_S \\ -(\bar{\boldsymbol{\Xi}}_T^d)^T \boldsymbol{Q}_T \end{bmatrix}. \tag{A8}$$

Converting Eq. (A8) into the state-space form, the following reduced model in the state-space form is obtained:

$$\begin{bmatrix} -\boldsymbol{\Lambda}_S^d & \boldsymbol{0} & \boldsymbol{0} \\ \boldsymbol{0} & \boldsymbol{I}_S^d & \boldsymbol{0} \\ \boldsymbol{0} & \boldsymbol{0} & -\boldsymbol{I}_T^d \end{bmatrix} \begin{bmatrix} \dot{\boldsymbol{q}}_S^d \\ \ddot{\boldsymbol{q}}_S^d \\ \dot{\bar{\boldsymbol{r}}}_T^d \end{bmatrix} + \begin{bmatrix} \boldsymbol{0} & \boldsymbol{\Lambda}_S^d & \boldsymbol{0} \\ \boldsymbol{\Lambda}_S^d & \boldsymbol{0} & -(\boldsymbol{\Phi}_S^d)^T \boldsymbol{K}_{ST}\bar{\boldsymbol{\Xi}}_T^d \\ \boldsymbol{0} & -(\bar{\boldsymbol{\Xi}}_T^d)^T \boldsymbol{K}_{TS}\boldsymbol{\Phi}_S^d & -\bar{\boldsymbol{\Theta}}_T^d \end{bmatrix} \begin{bmatrix} \boldsymbol{q}_S^d \\ \dot{\boldsymbol{q}}_S^d \\ \bar{\boldsymbol{r}}_T^d \end{bmatrix} = \begin{bmatrix} \boldsymbol{0} \\ (\boldsymbol{\Phi}_S^d)^T \boldsymbol{f}_S \\ -(\bar{\boldsymbol{\Xi}}_T^d)^T \boldsymbol{Q}_T \end{bmatrix}. \tag{A9}$$

Consequently, the proposed formulation in Section 4 can be also obtained from the second-order governing equation.

| Selected model | Computational time [sec] |
|:---:|:---:|
| **Mode superposition model** | 541.808921 |
| **Uncoupled model** | 80.458196 |
| **Proposed model** | 95.629184 |

Table. 1. Computational costs for constructing reduced models in the 3d flange-pipe problem.

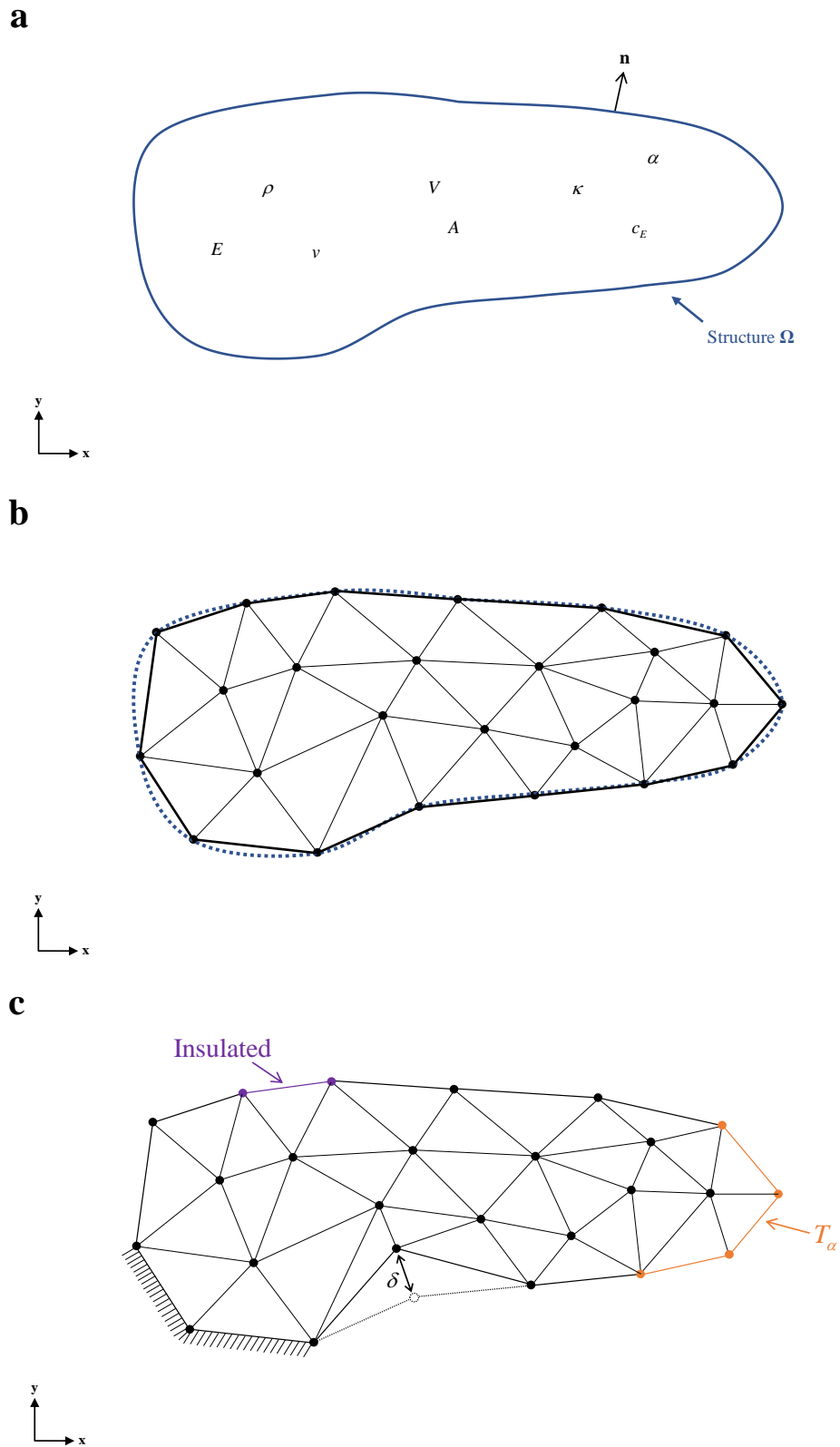

Fig. 1. Fully coupled thermoelastic problem for an example structure $\Omega$ using the FE discretization : (a) shape and material properties of the structure $\Omega$, (b) a FE discretization of the structure $\Omega$ and (c) boundary conditions of thermomechanical problem in the structure $\Omega$

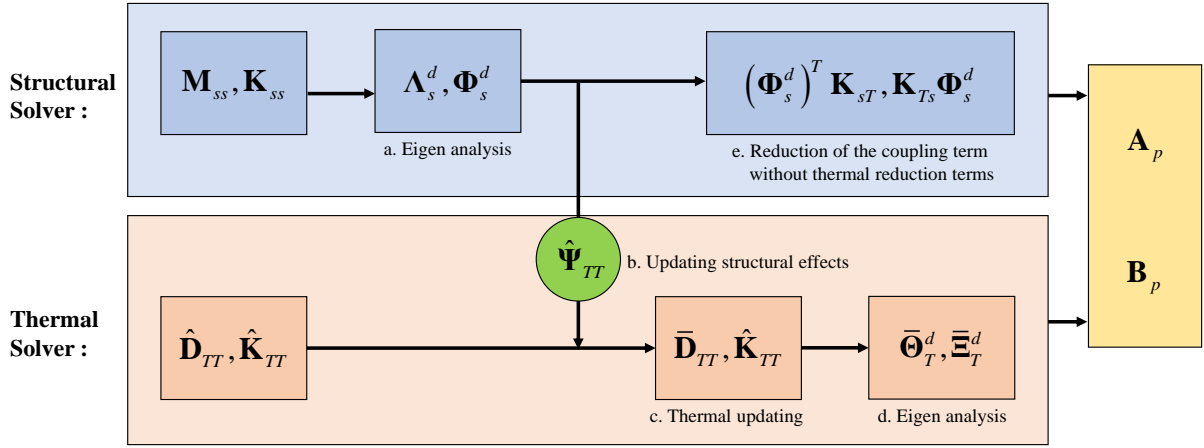

Fig. 2. The process of the constructing the reduced model with the proposed method

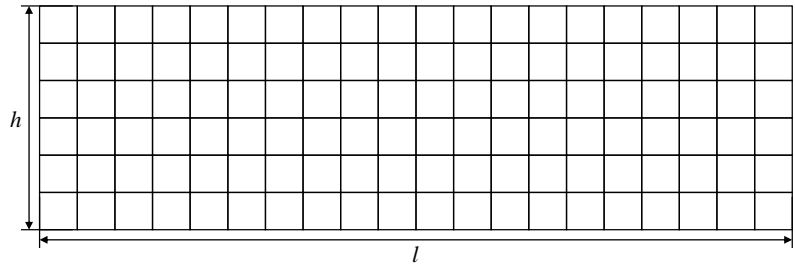

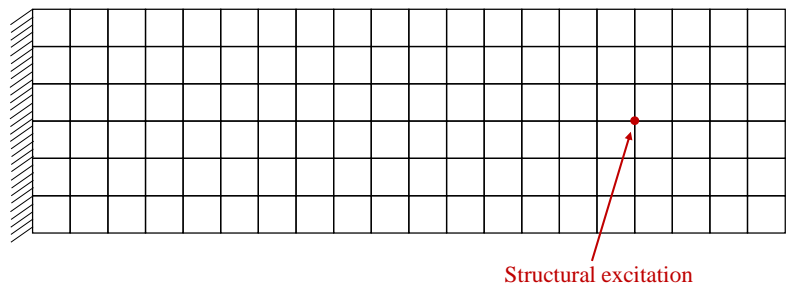

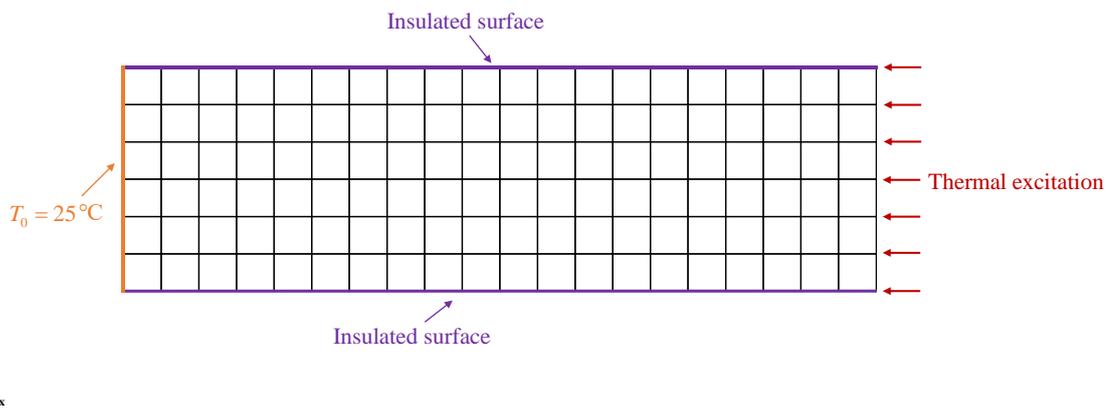

Fig. 3. Geometry, boundary and excitation conditions for the 2d plate problem : (a) a geometry condition of the 2d plate, (b) boundary conditions and an excitation node for the structural part in the 2d plate problem and (c) boundary conditions and an excitation nodes for the thermal part in the 2d plate problem.

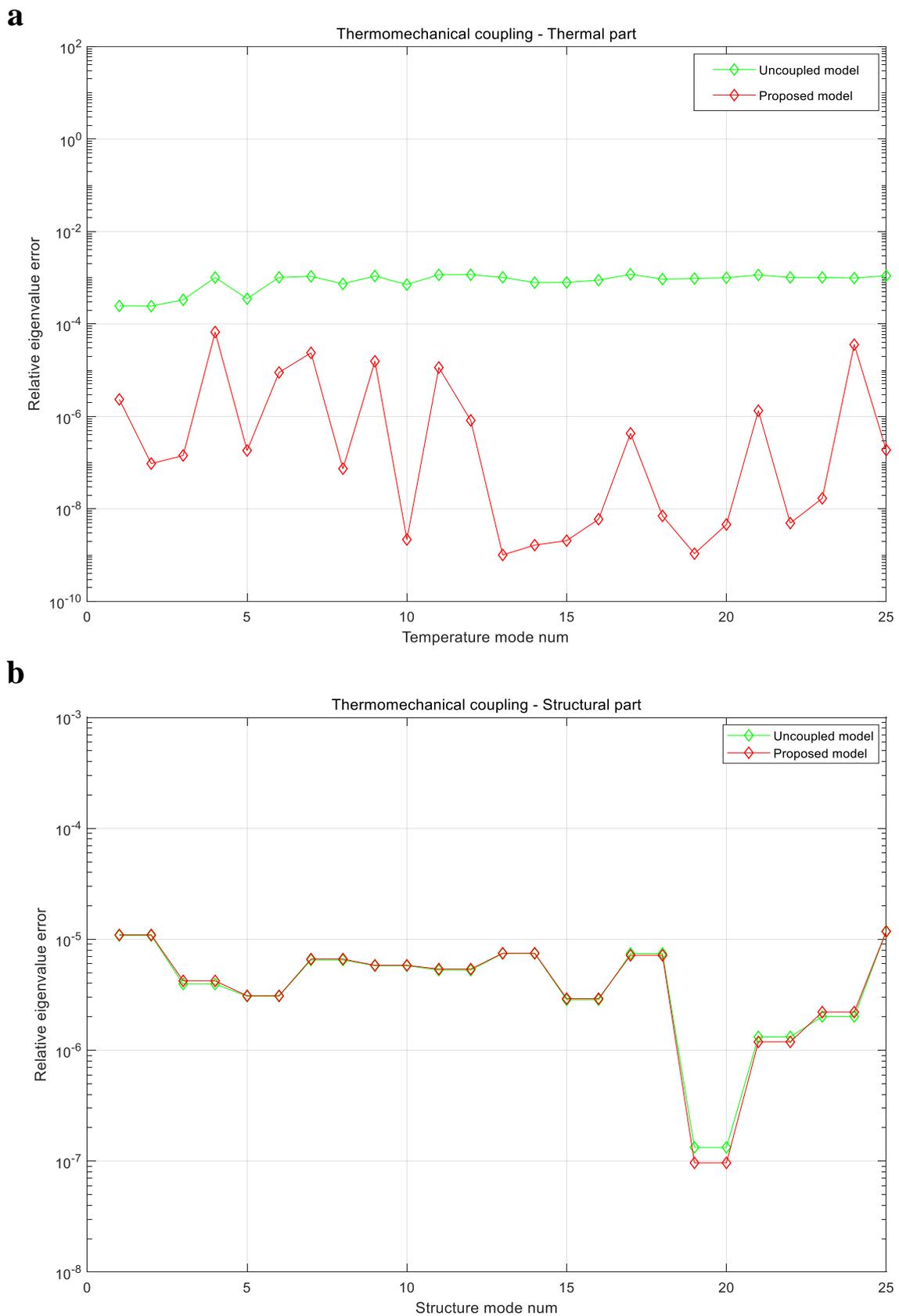

Fig. 4. Results of the relative eigenvalue errors of the 2d plate problem : (a) relative eigenvalue errors of the thermal part and (b) relative eigenvalue errors of the structural part.

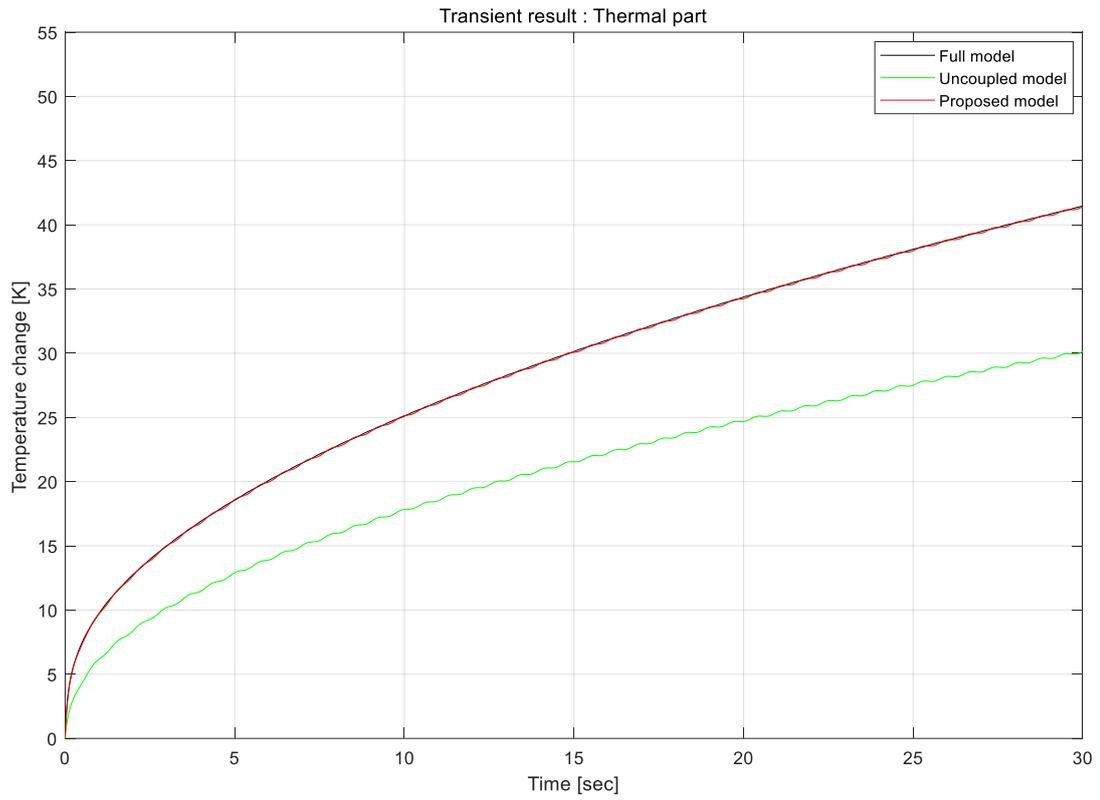

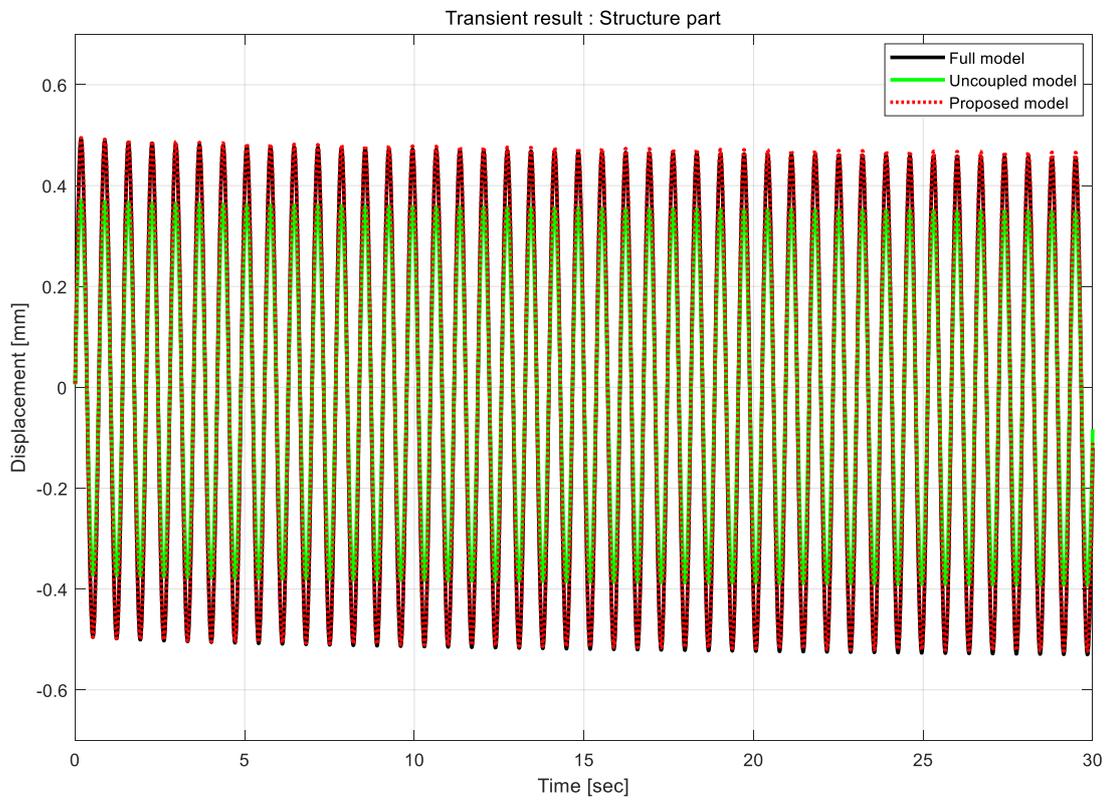

Fig. 5. Transient results for the 2d plate problem : (a) maximum temperature change results and (b) maximum displacement results.

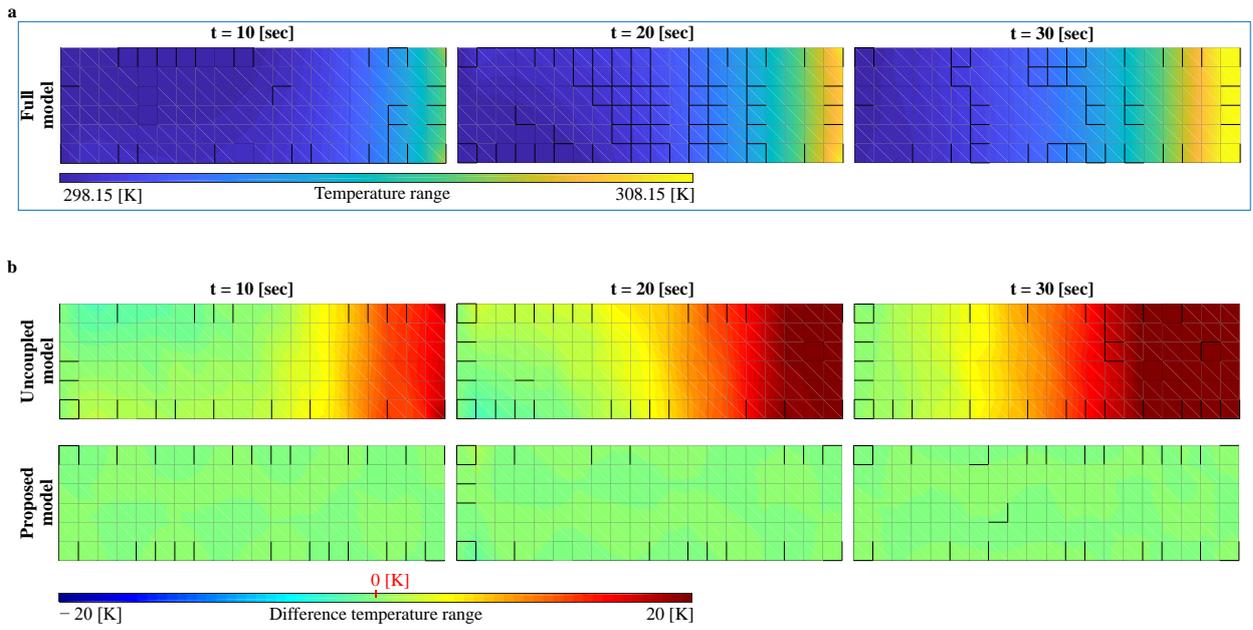

Fig. 6. Temperature fields of the 2d plate problem : (a) temperature fields of the full model in the 2d plate problem and (b) differences of the temperature fields between the full and each reduced model in the 2d plate problem.

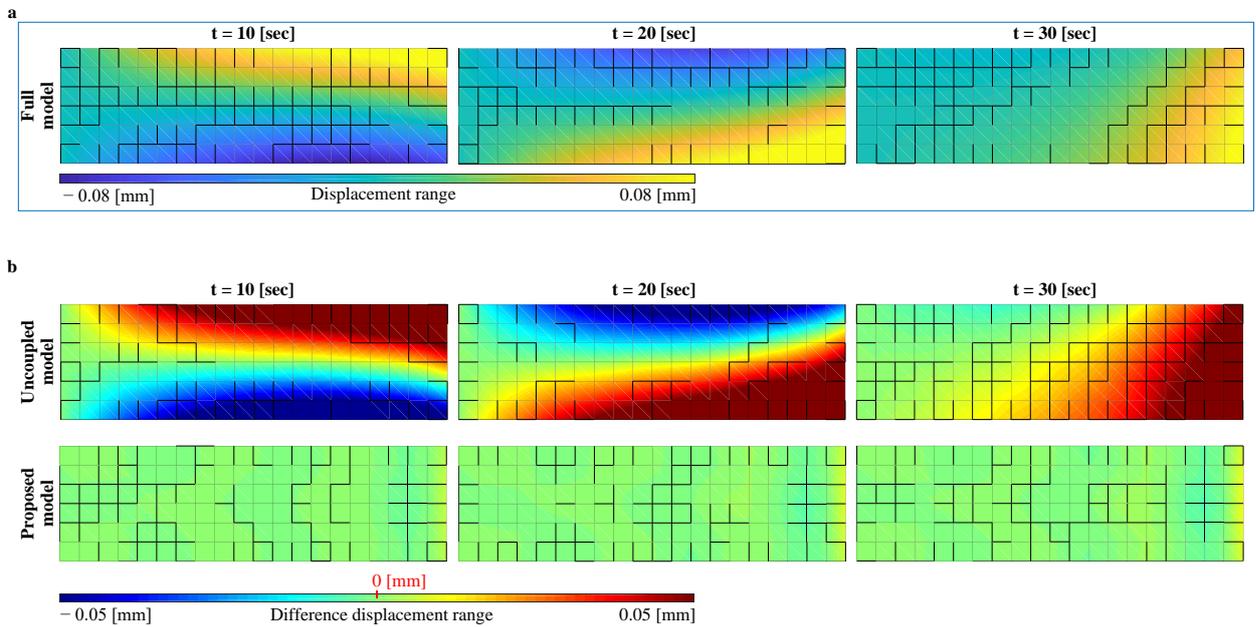

Fig. 7. x direction displacement fields of the 2d plate problem : (a) x direction displacement fields of the full model in the 2d plate problem and (b) differences of the x direction displacement fields between the full and each reduced model in the 2d plate problem.

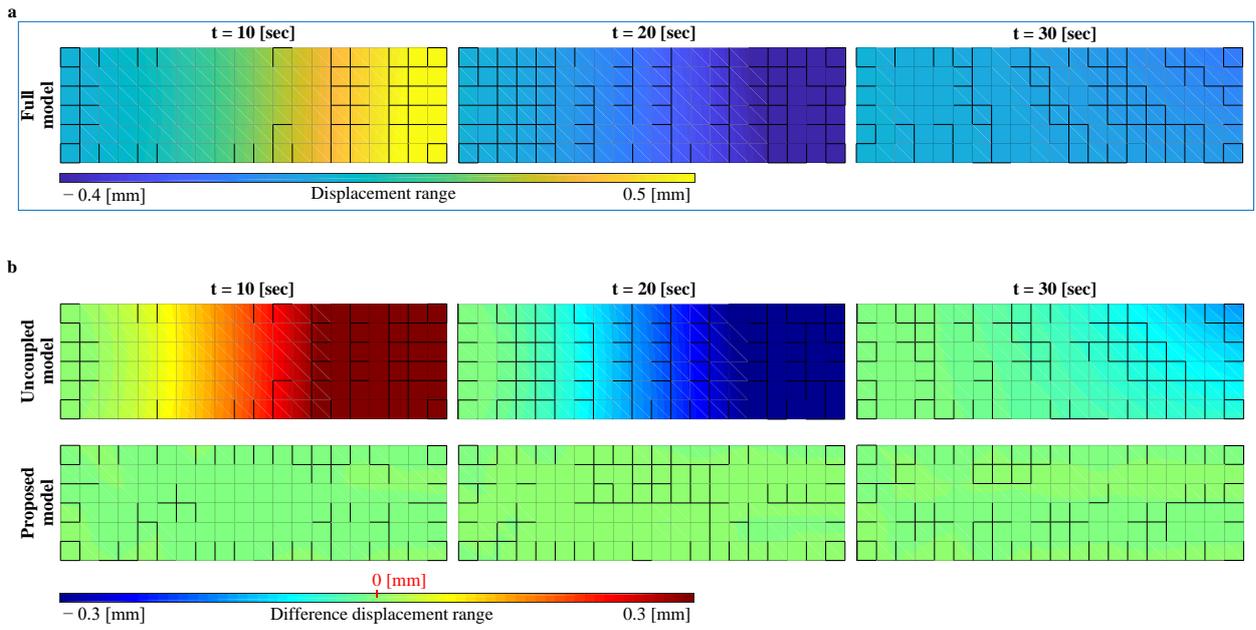

Fig. 8. y direction displacement fields of the 2d plate problem : (a) y direction displacement fields of the full model in the 2d plate problem and (b) differences of the y direction displacement fields between the full and each reduced model in the 2d plate problem.

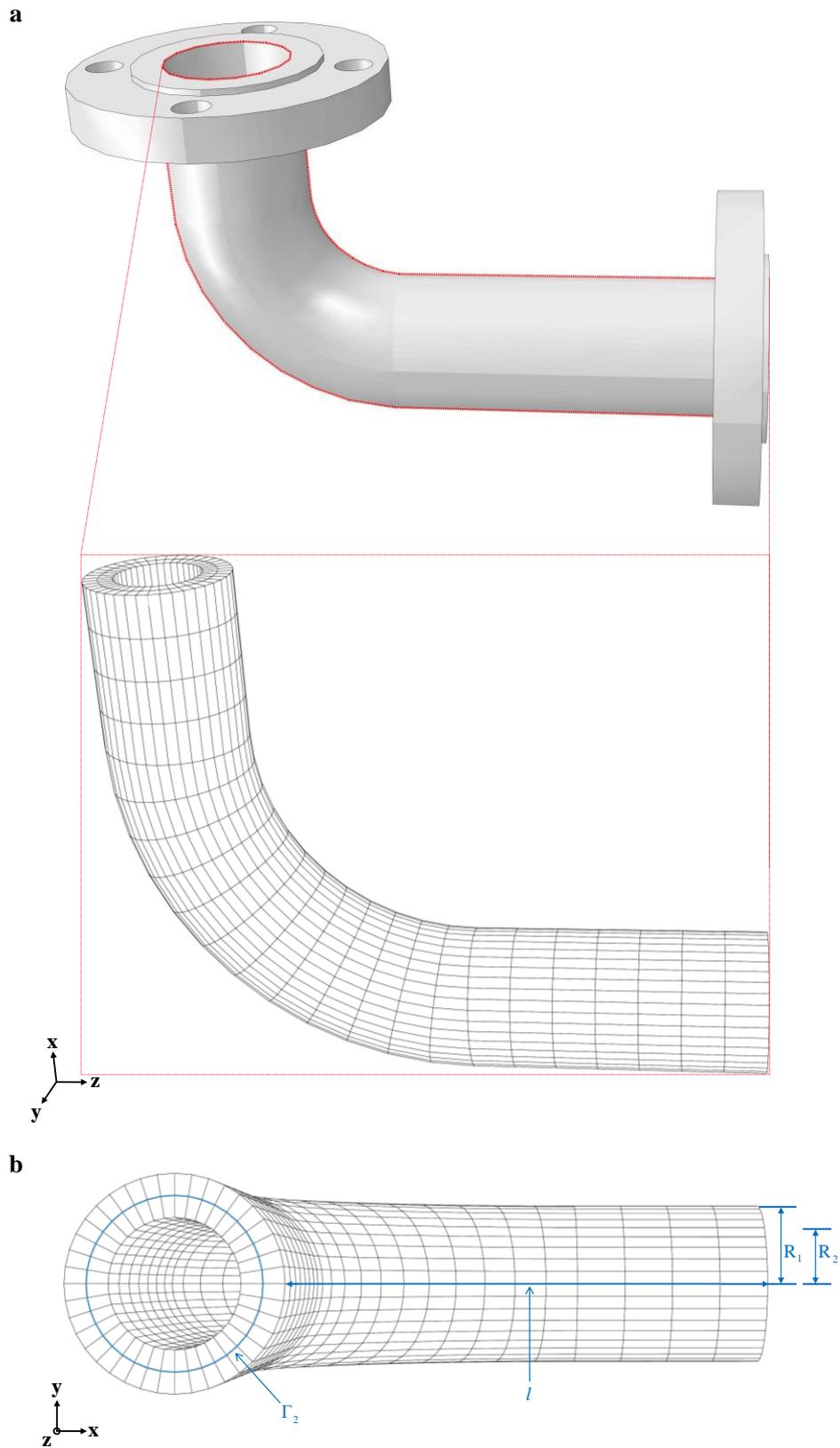

Fig. 9. Shape, FE modeling result and geometry quantities for the flange-pipe problem : (a) a shape of the 3d flange-pipe problem and its FE modeling results and (b) geometry quantities ($R_1$, $R_2$ and $l$) of the 3d flange-pipe problem and a location of a surface $\Gamma_2$.

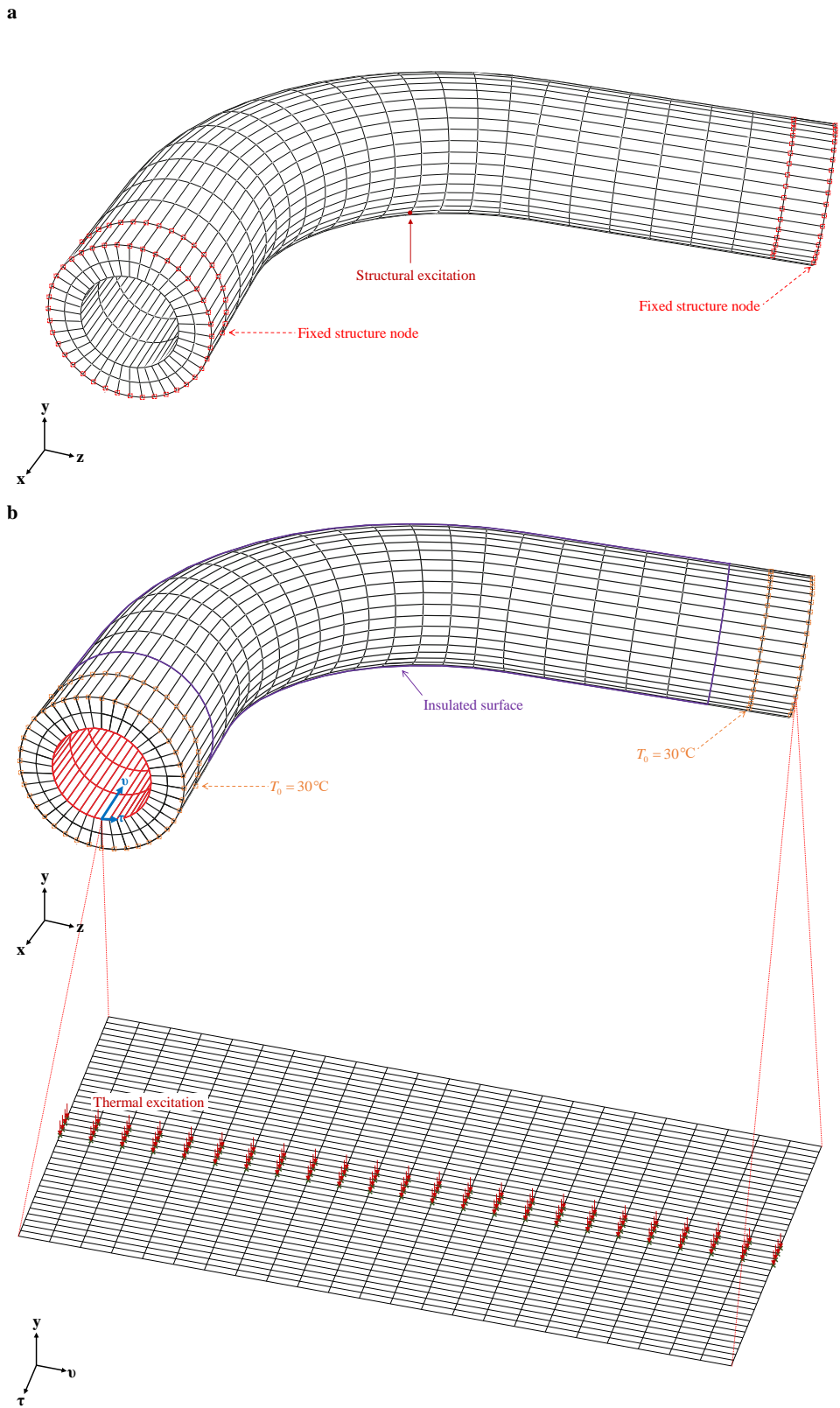

Fig. 10. Boundary conditions and excitation locations for the 3d flange-pipe problem : (a) boundary conditions, remained structure nodes and an excitation condition for the structure part in the 3d flange-pipe problem and (b) boundary conditions, remained temperature nodes and an excitation condition for the thermal part in the 3d flange-pipe problem.

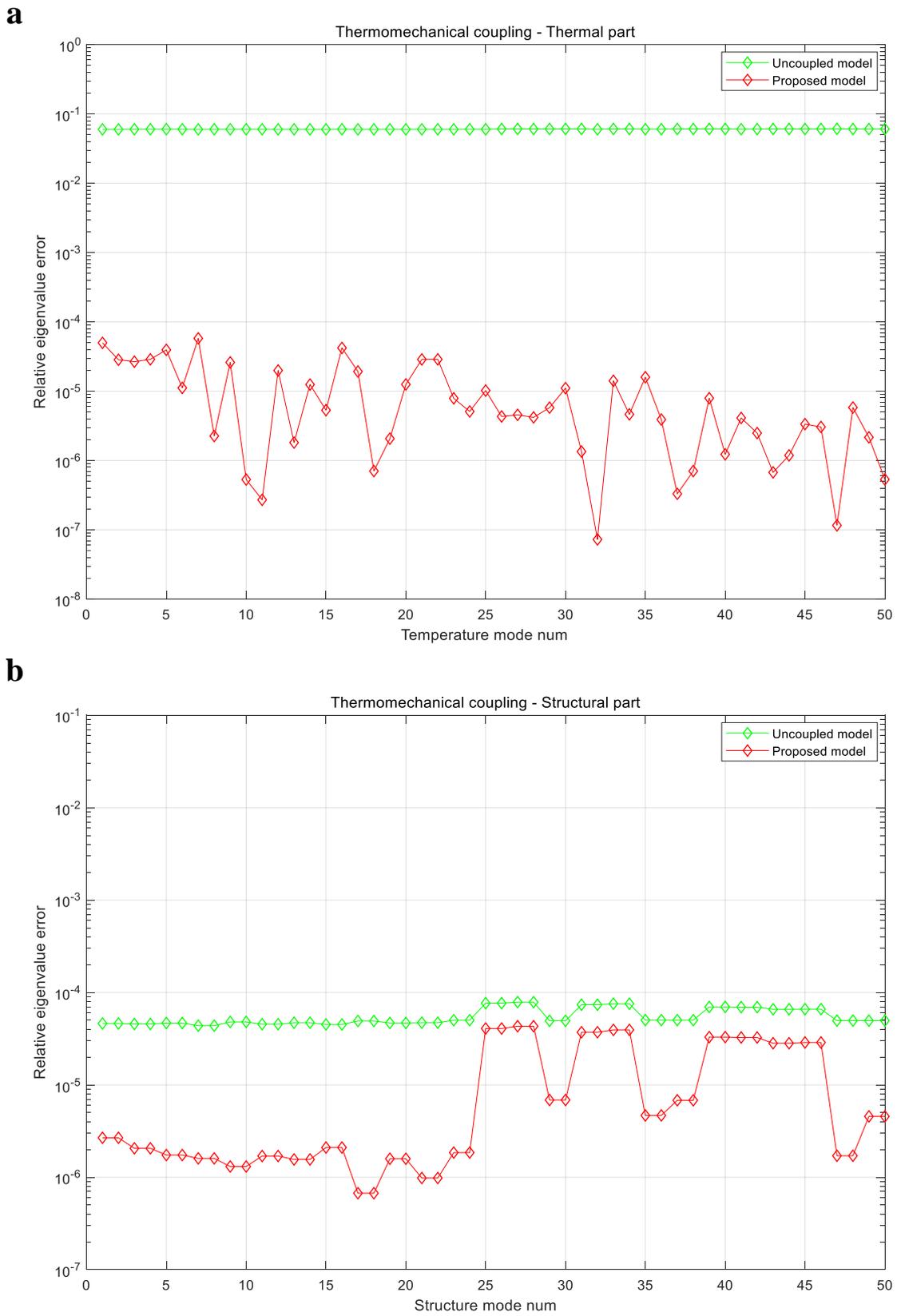

Fig. 11. Results of the relative eigenvalue errors of the 3d flange-pipe problem: (a) relative eigenvalue errors of the thermal part and (b) relative eigenvalue errors of the structure part.

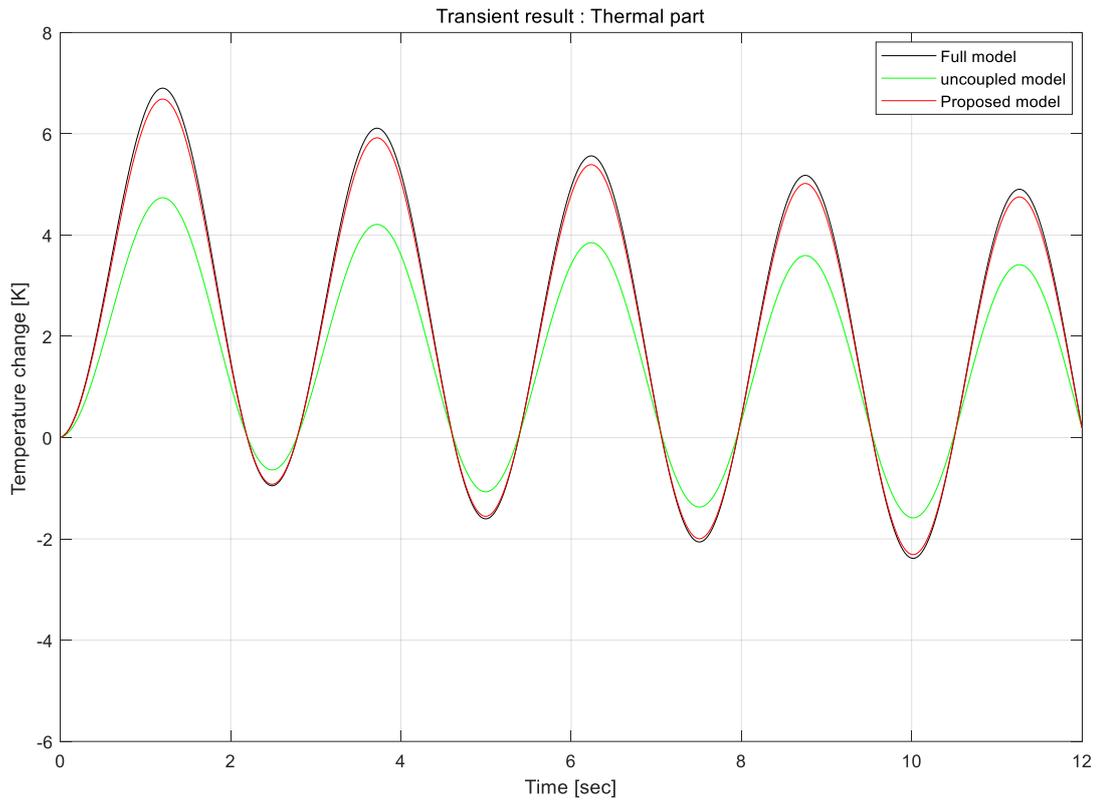

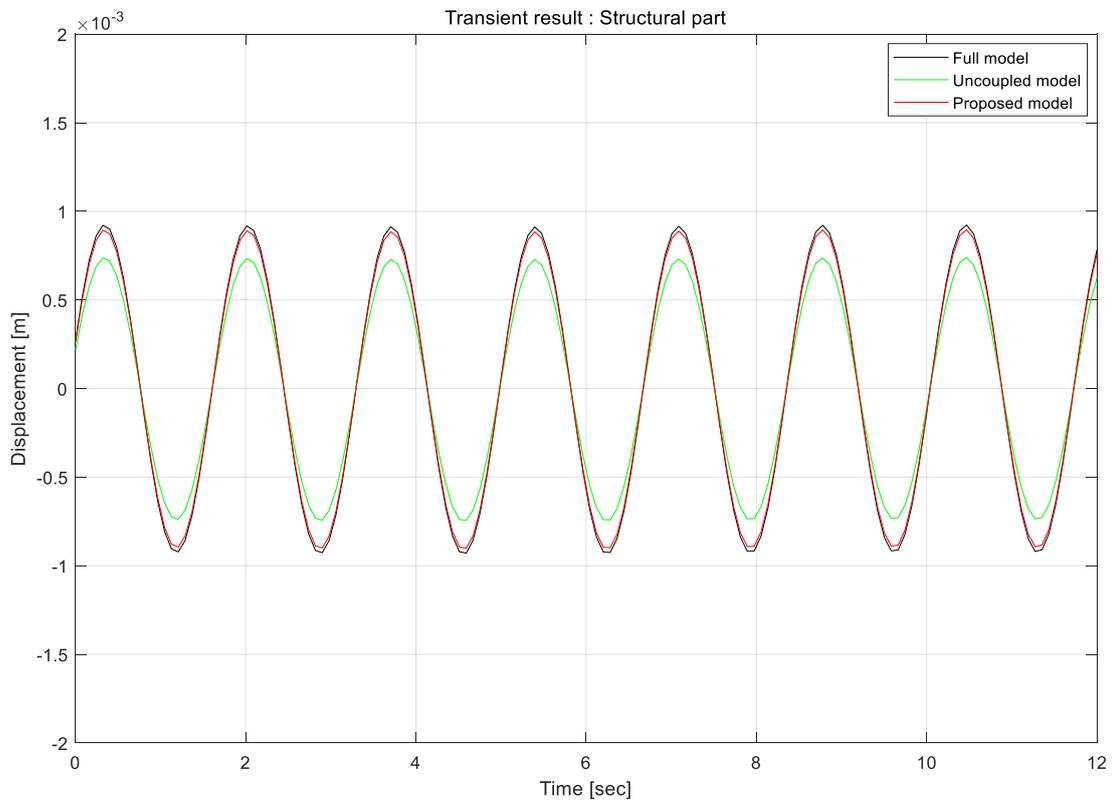

Fig. 12. Transient results for the 3d flange-pipe problem : (a) maximum temperature change results and (b) maximum displacement results.

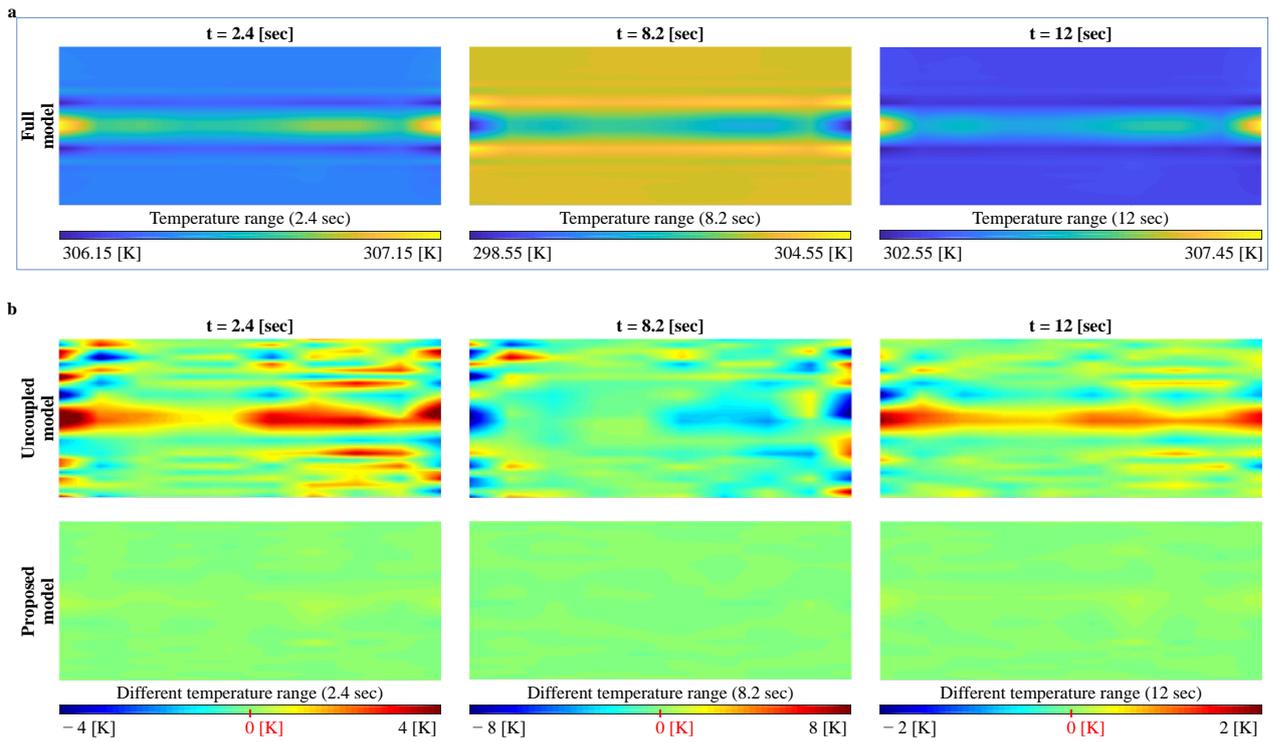

Fig. 13. Temperature fields of the 3d flange pipe problem in the surface $\Gamma_2$ : (a) temperature fields of the full model and (b) differences of the temperature fields between the full and each reduced model.

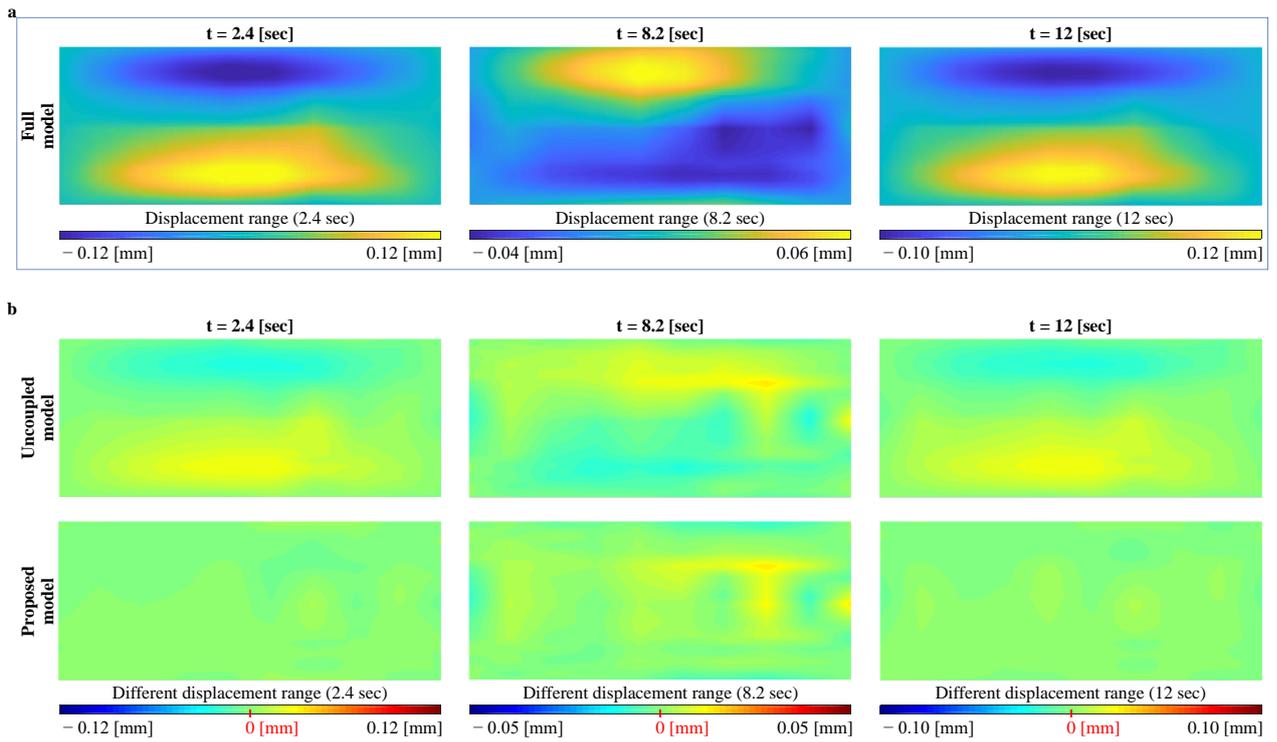

Fig. 14. X direction displacement fields of the 3d flange pipe problem in the surface $\Gamma_2$ : (a) x direction displacement fields of the full model and (b) differences of x direction displacement fields between the full and reduced model

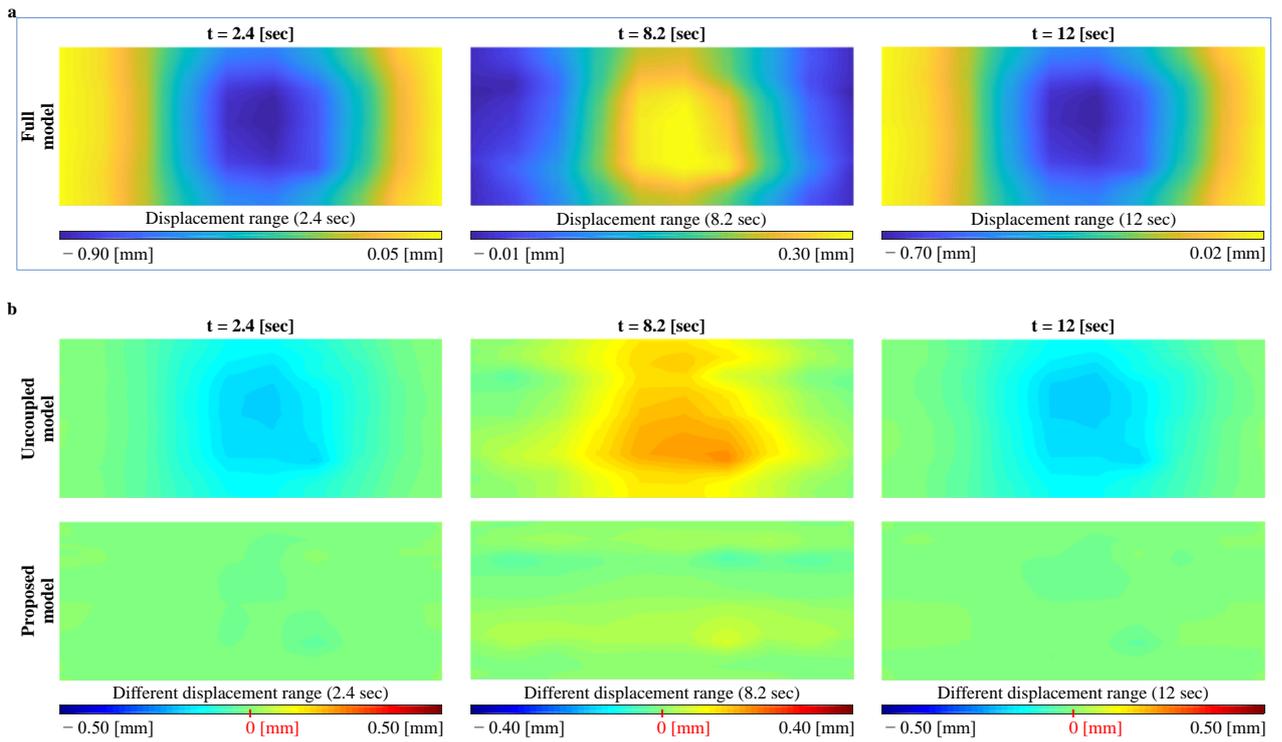

Fig. 15. Y direction displacement fields of the 3d flange pipe problem in the surface $\Gamma_2$ : (a) y direction displacement fields of the full model and (b) differences of y direction displacement fields between the full and each reduced model

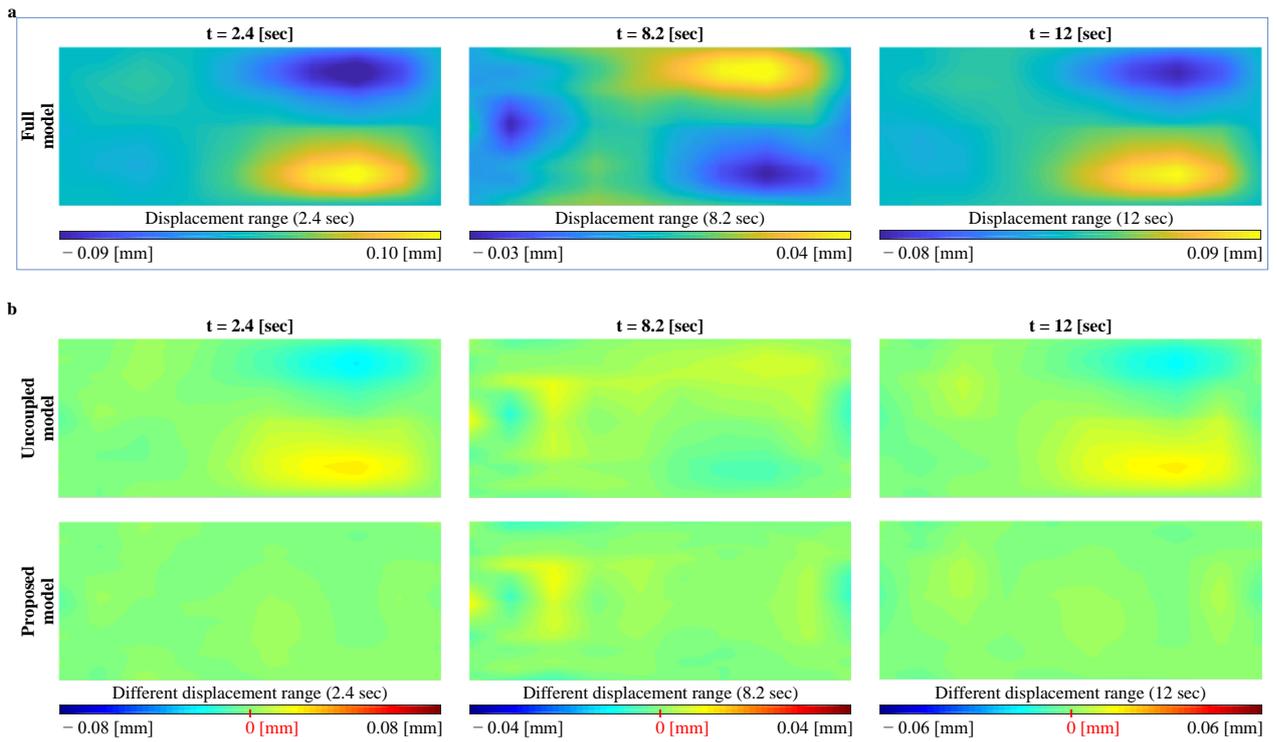

Fig. 16. Z direction displacement fields of the 3d flange-pipe problem in the surface $\Gamma_2$ : (a) z direction displacement fields of the full model and (b) differences of z direction displacement fields between the full and each reduced model.

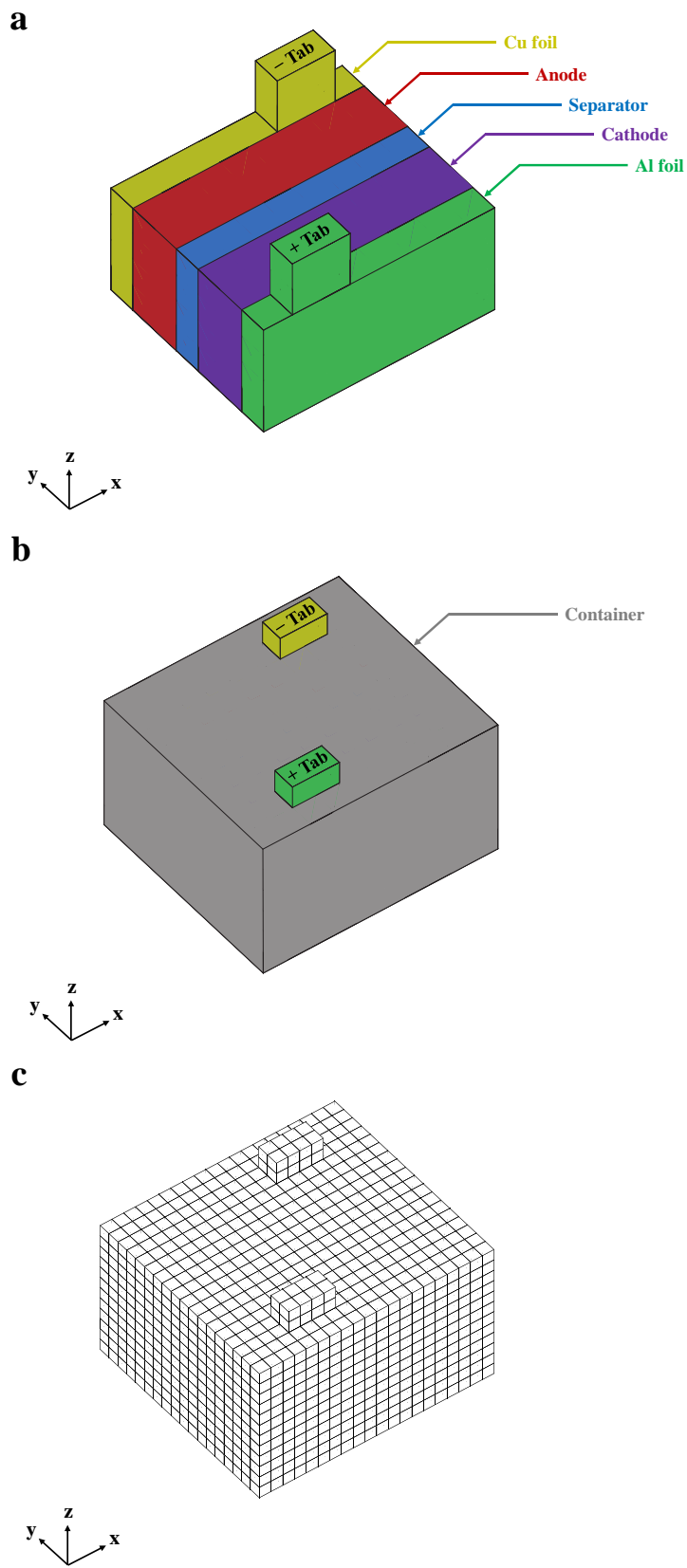

Fig. 17. A structure of the lithium-ion battery and its mesh result : (a) a structure of the pouch cell in the lithium-ion battery, (b) an entire structure of the lithium-ion battery and (c) a mesh result of the lithium-ion battery problem.

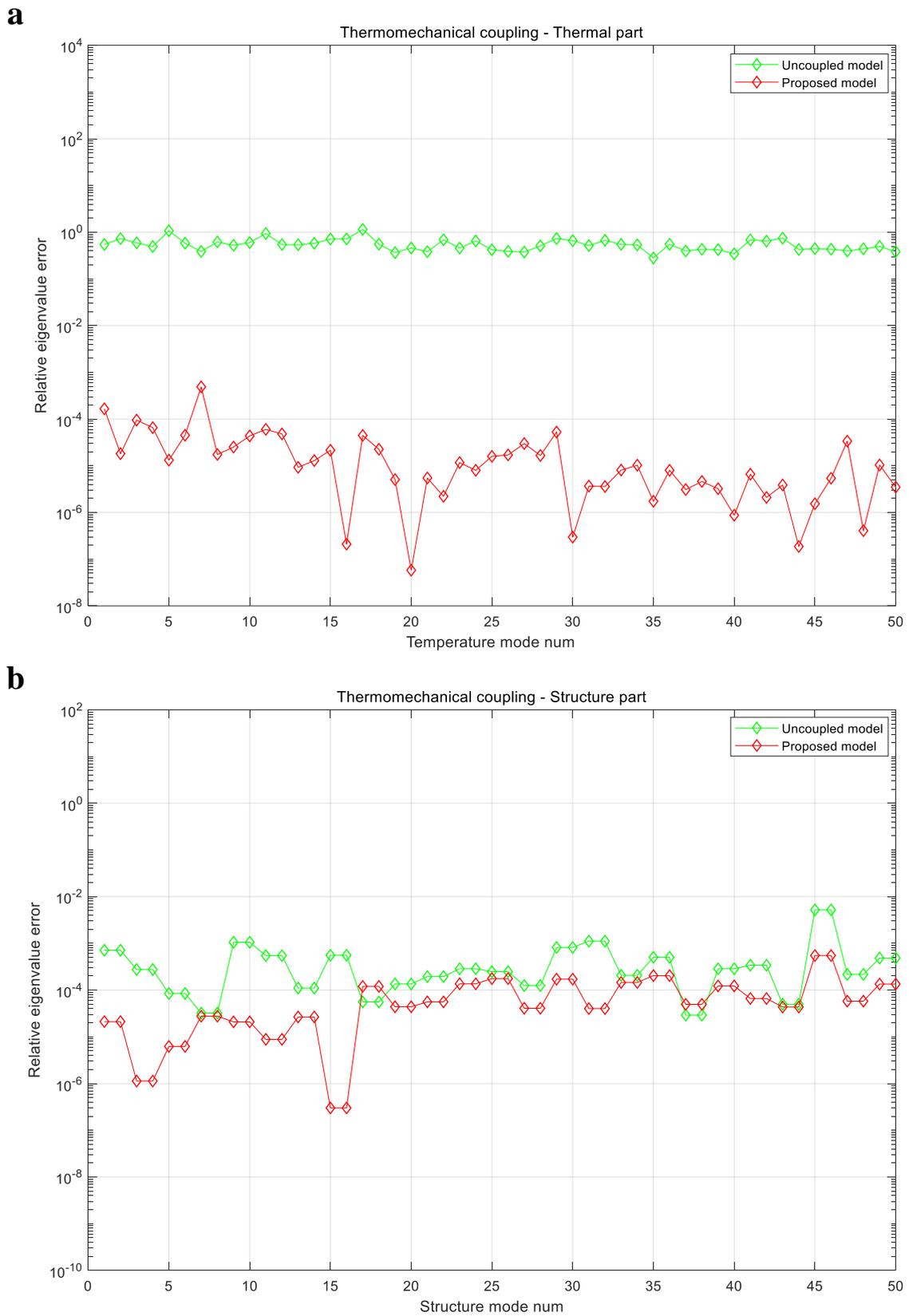

Fig. 18. Results of the relative eigenvalue errors of the lithium-ion battery problem: (a) relative eigenvalue errors of the temperature part and (b) relative eigenvalue errors of the structure part.

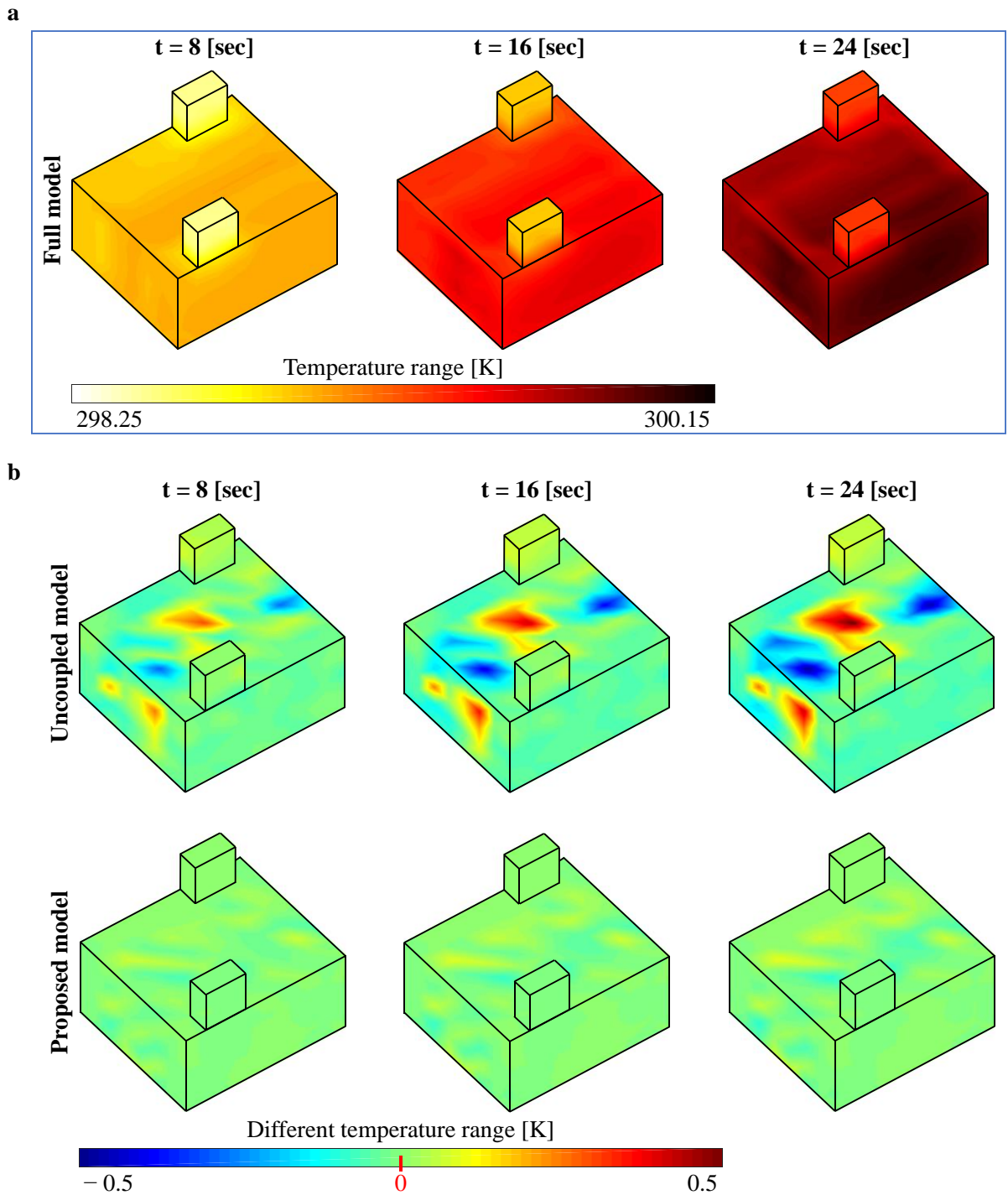

Fig. 19. Temperature fields of the pouch cell in the lithium-ion battery : (a) temperature fields of the full model and (b) differences of the temperature fields between the full and each reduced model.

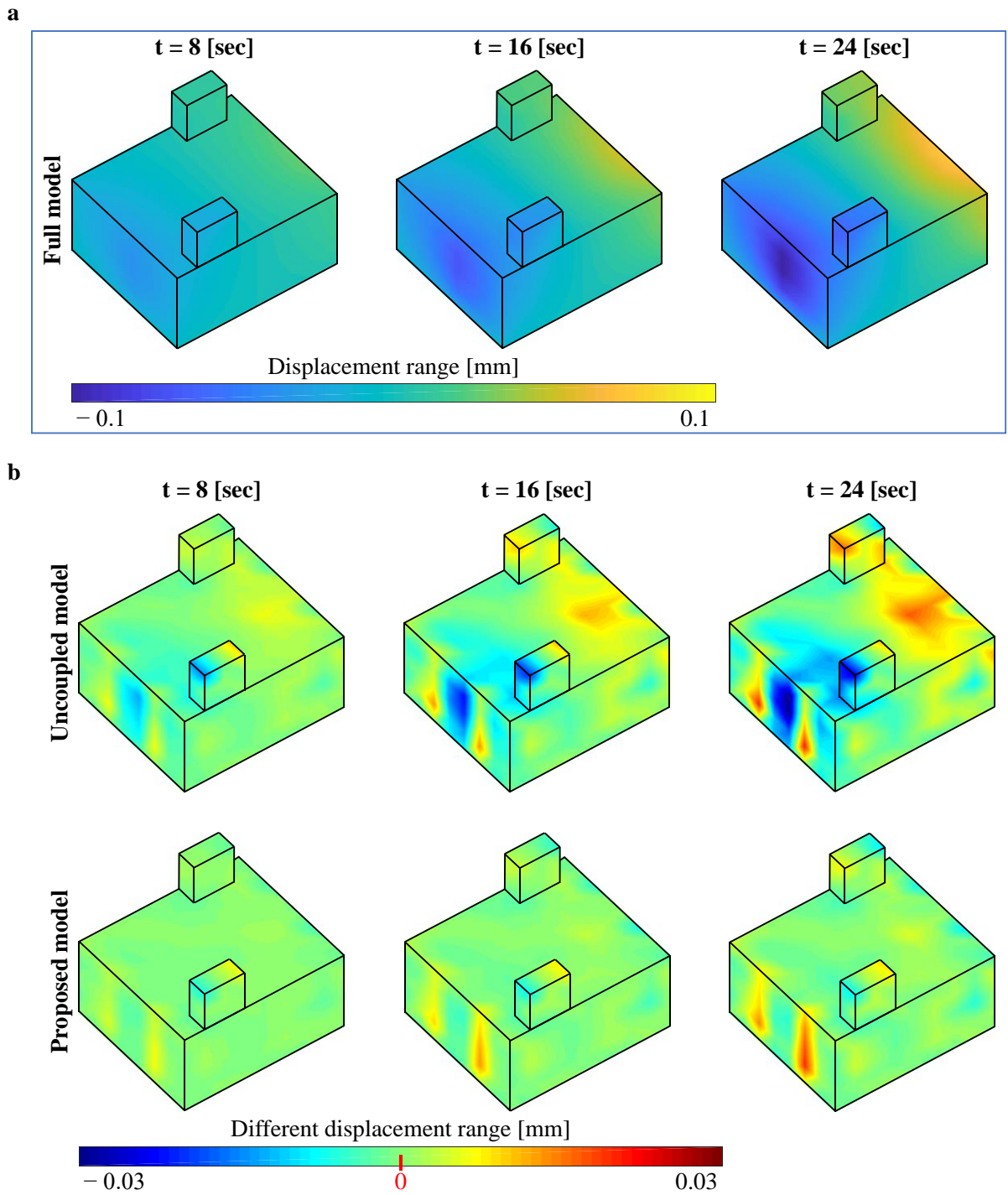

Fig. 20. X direction displacement fields of the pouch cell in the lithium-ion battery : (a) displacement fields of the full model and (b) differences of the displacement fields between the full and each reduced model.

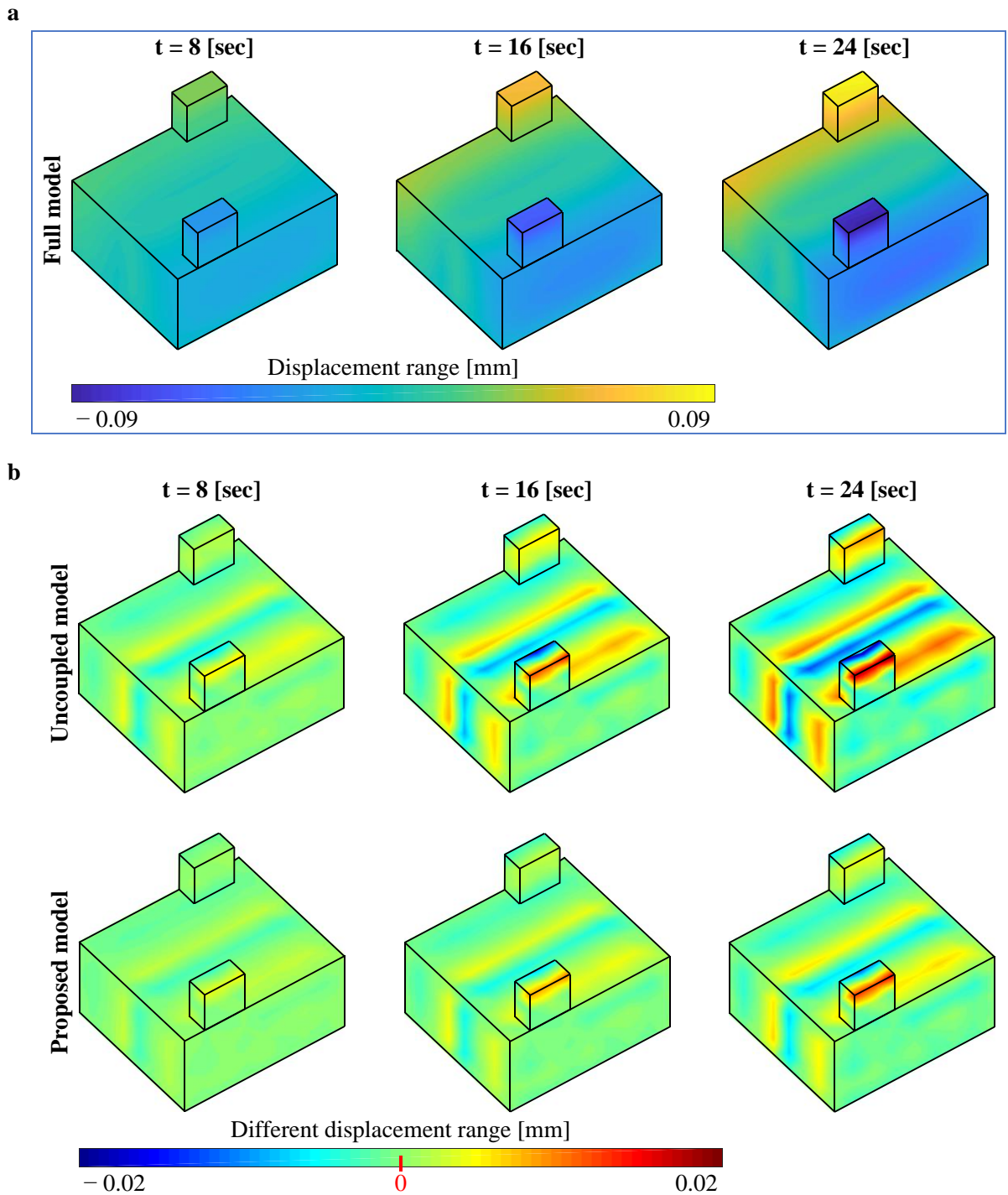

Fig. 21. Y direction displacement fields of the pouch cell in the lithium-ion battery : (a) displacement fields of the full model and (b) differences of the displacement fields between the full and each reduced model.

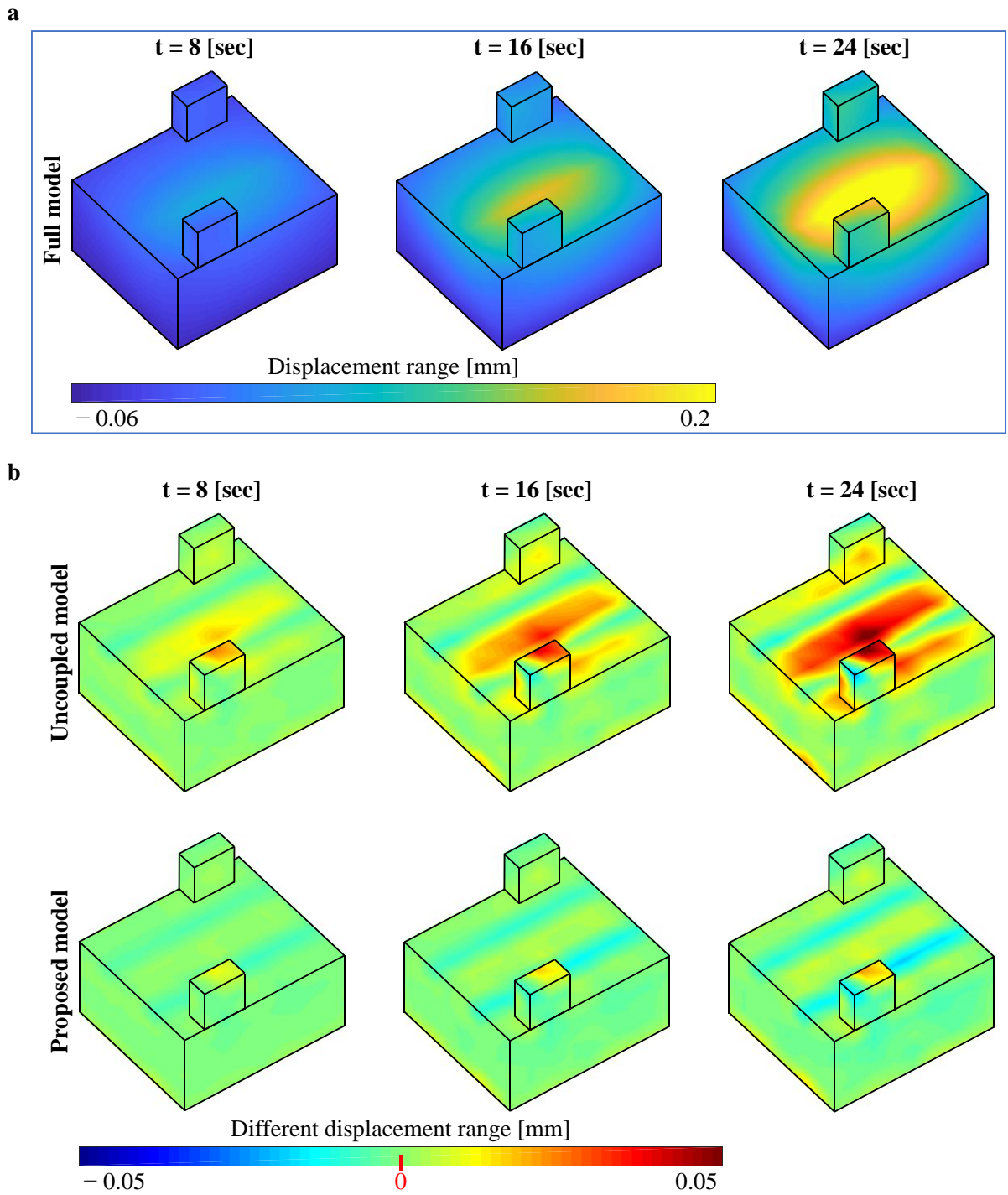

Fig. 22. Z direction displacement fields of the pouch cell in the lithium-ion battery : (a) displacement fields of the full model and (b) differences of the displacement fields between the full and each reduced model.

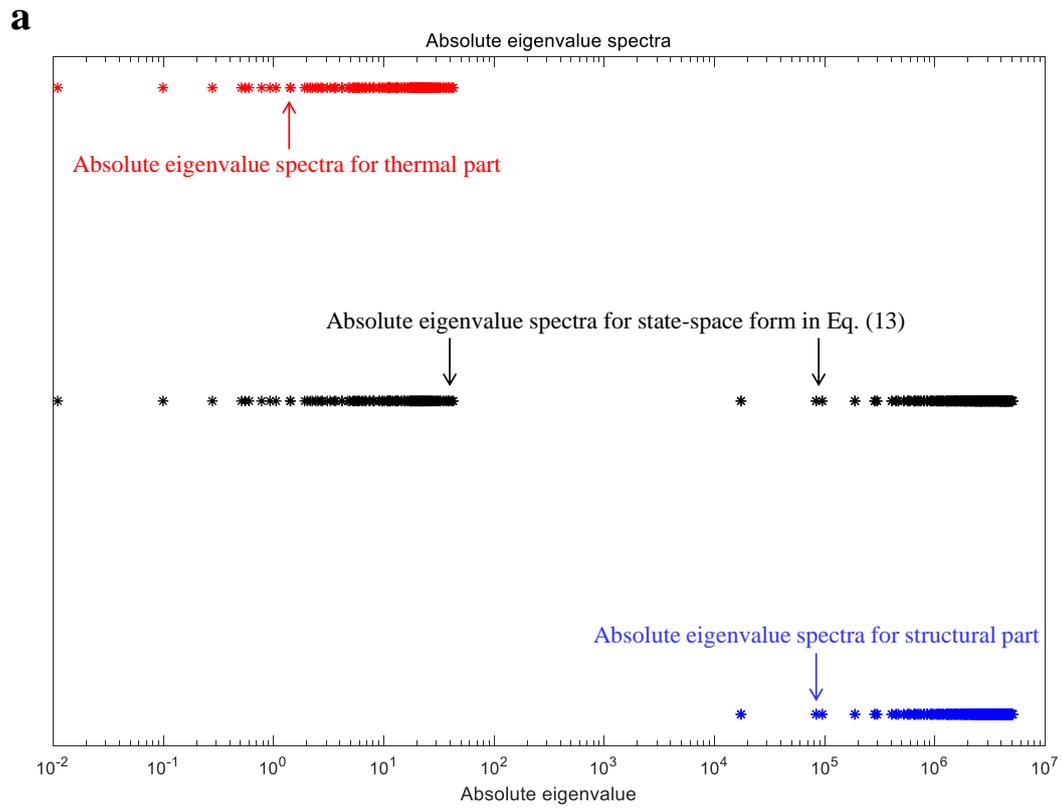

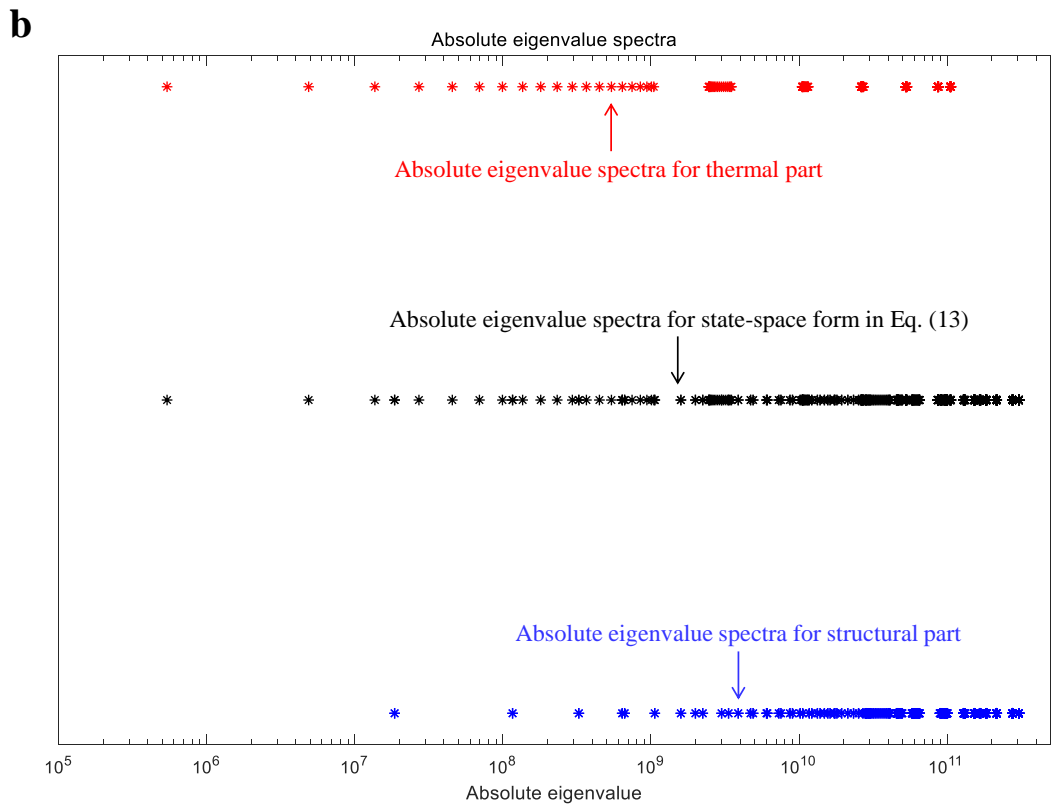

Fig. 23. Results of the absolute eigenvalue spectra in the 2d plate problem: (a) absolute eigenvalue spectra at macro-scale ($h = 0.042m$, $l = 0.140m$) and (b) absolute eigenvalue spectra at micro-scale ($h = 4\mu$m, $l = 20\mu$m)

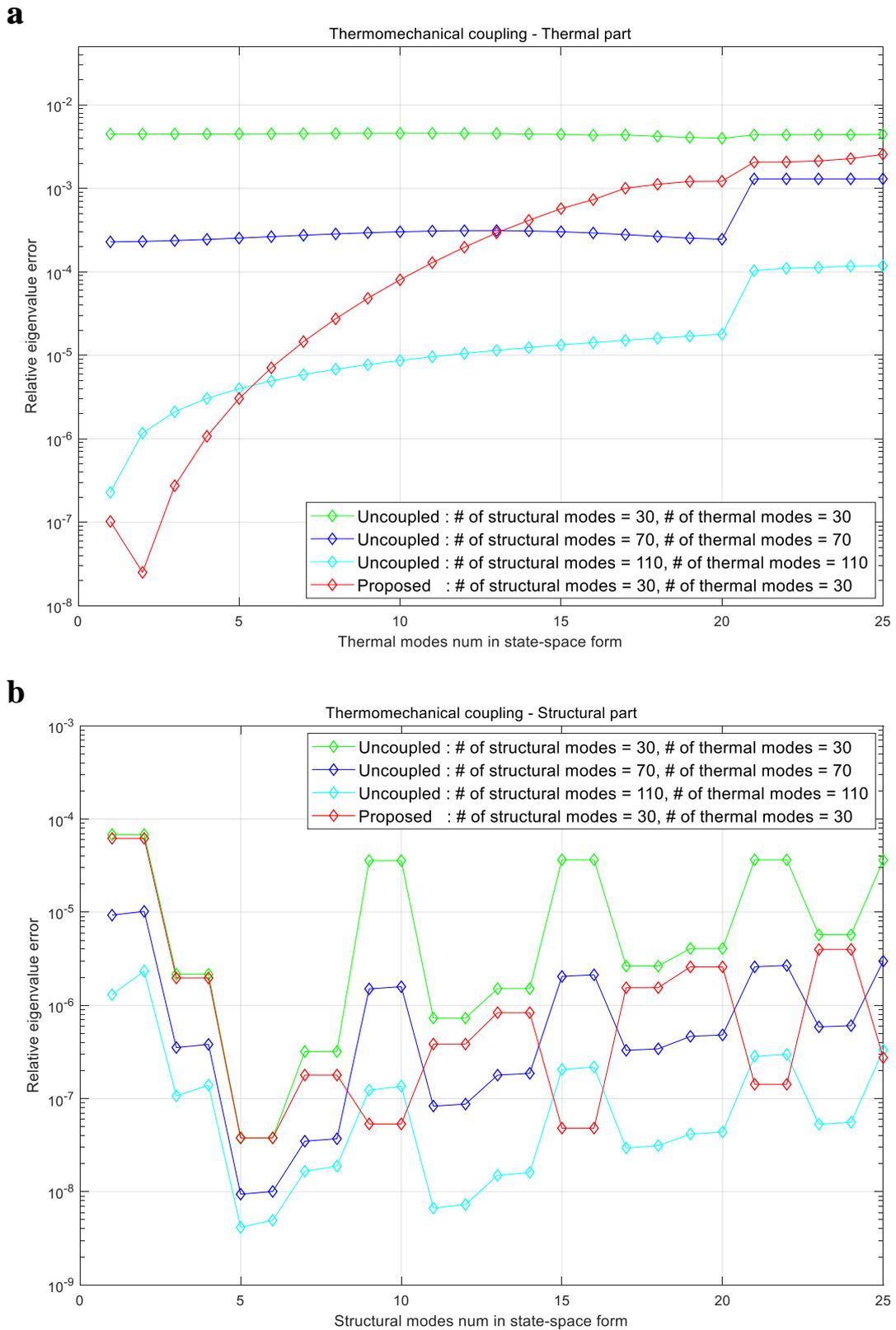

Fig. 24. Results of the relative eigenvalue error of the 2d plate problem at micro-scale ($h = 4\mu m$, $l = 20\mu m$) : (a) relative eigenvalue errors of the thermal part and (b) relative eigenvalue errors of the structural part. (For interpretation of the references to colors in this figure legend, the reader is referred to the web version of this article)